\documentclass{amsart}
\usepackage{latexsym}
\usepackage{amstext}
\usepackage{amsmath}
\usepackage{amssymb}
\usepackage{amsthm}
\usepackage{enumitem}
\usepackage{mathrsfs} 
\usepackage{url}
\usepackage{bm}
\usepackage{comment}
\usepackage[utf8]{inputenc}
\usepackage[T1]{fontenc}
\usepackage{amsmath,amssymb,amsthm,mathtools}
\usepackage[colorlinks={true},linkcolor={red},citecolor={blue}, urlcolor={blue}]{hyperref}
\usepackage[margin=2.5 cm]{geometry}

\theoremstyle{plain}
\newtheorem{theorem}{Theorem}[section]

\newtheorem{corollary}[theorem]{Corollary}
\newtheorem{lemma}[theorem]{Lemma}

\theoremstyle{definition}
\newtheorem{definition}[theorem]{Definition}

\newtheorem{remark}[theorem]{Remark}

\newtheorem{example}[theorem]{Example}
\newtheorem{proposition}[theorem]{Proposition}

\newcommand{\K}{\mathbb{K}}

\newcommand{\X}{\mathbb{X}}

\allowdisplaybreaks
\title[Roth-Type Diophantine Inequalities over Adelic Curves]{Explicit Bounds for the Number of Solutions to Roth-Type Inequalities over Adelic Curves}
\author{YuXi Li, GuanHeng Zhao}
\begin{document}
	\begin{abstract}
		In this paper, we provide an upper bound on the number of a special system of Roth-type Diophantine
		inequalities over proper adelic curves. By refining the technique from previous work, we successfully eliminate
		the large-height restriction on the solutions. Furthermore, the constraints on the parameters within the
		inequality are further relaxed. We also give a constructive quantitative generalization of the Vojta-Dolce-Zucconi system, replacing compactness-based extraction by a discrete pigeonhole argument and a finite grid partition, so that the resulting bound has explicit dependence on the relevant parameters.
	\end{abstract}	

	\maketitle	
	\tableofcontents
	
	\section{Introduction}\label{cha_intro}
	\subsection{Background}
	An \emph{adelic curve} in the sense of Chen--Moriwaki \cite{Ch-Mo} is a triple
	$$
	\X=(\K,\Omega,\mu),
	$$
	where $\K$ is a field, $\Omega$ is a measure space parametrizing a family of absolute values $\omega\mapsto|\cdot|_\omega$ on $\K$, and $\mu$ is a measure on $\Omega$. Global fields arise canonically as adelic curves; more generally, any countable field admits non-trivial adelic curve structures. Thus an adelic curve may be viewed as a field equipped with a measure-parametrized family of absolute values. This notion naturally unifies number fields, function fields, and arithmetic function fields---the latter being finitely generated extensions of $\mathbb{Q}$ endowed with absolute values derived from Arakelov-completed arithmetic varieties (see \cite[Section 3]{vo}). This unifying framework motivates the study of Diophantine problems over adelic curves, since it allows theorems from traditionally distinct settings to be proved in a uniform way.
	
	The Diophantine approximation problem over $\X$ concerns, for fixed $\alpha_1,\dots,\alpha_q\in\K$, the quantities
	$$
	-\log^-|\xi-\alpha_j|_\omega,
	$$
	their integrals over $\Omega$, and the comparison of these integrals with the logarithmic height $h_\X(\xi)$. Classically, Roth's theorem says, roughly, that for a fixed algebraic $\alpha$ and every $\varepsilon>0$, the inequality
	$$
	-\log|\xi-\alpha|>(2+\varepsilon)h(\xi)
	$$
	has only finitely many solutions $\xi\in\overline{\mathbb{Q}}$. 
	
	Quantitative versions of Roth's theorem are known for number fields and function fields \cite{DR,bom,bom1}; see also \cite{gr1,sc2,C} for related quantitative results. For adelic curves, Dolce and Zucconi generalized Roth's theorem to proper adelic curves under technical conditions on logarithmic absolute values \cite{dz}, recovering as special cases Corvaja's results for fields with a product formula and Vojta's arithmetic function field results. More recently, Dolce studied quantitative solutions of certain Diophantine inequalities over adelic curves \cite{d}, providing a method for handling quantitative versions of Roth's theorem in this setting. By modifying and extending the methods introduced in \cite{d}, we substantially strengthen the result of \cite{d}. We first fix the notation that will be used throughout the paper; further details concerning these definitions are provided in Section \ref{cha:prel}.
	\subsection{Main results}
	\textbf{Notation.} Throughout this paper, we adopt the following conventions:
	\begin{itemize}
		\item $\mathbb{K}$ is a field of characteristic $0$, and $\mathbb{X}=(\mathbb{K},\Omega,\phi)$ is a proper adelic curve, meaning that it satisfies a product formula.
		\item $\Omega$ is a measure space equipped with a measure $\mu$, which parametrizes the absolute values of $\mathbb{K}$ via $\phi$.
		\item $h=h_{\mathbb{X}}$ denotes the logarithmic height on $\mathbb{X}$, and $H=e^h$ is the exponential height.
		\item $S\subset\Omega$ is a fixed subset of finite measure. For $i=1,\dots,n$, the sets $S_i$ are measurable and satisfy
		$$ S=\bigsqcup_{i=1}^{n} S_i. $$
		\item $\alpha_1,\ldots,\alpha_n\in\mathbb{K}^{\times}$ are fixed, distinct elements.
		\item $A>0$ and $\varepsilon>0$ are fixed real constants.
		\item For every $X>0$, define $c(X)$ by
		$$ c(X)=(2+\varepsilon)\log X-2\log 2. $$
		\item For every $\delta>0$, define $N(\delta)$ by
		$$ N(\delta)=\max\left(\lceil 21^2\log 2n\rceil, \left\lceil \frac{49\log 2n}{\delta^2}\right\rceil\right), $$
		and define $b(\delta)$ and $B(\delta)$ by
		$$ B(\delta)=4^{2N(\delta)!}\prod_{i=1}^n H(\alpha_i)^{\frac{2N(\delta)!}{n}},\qquad b(\delta)=\log B(\delta). $$
		\item $\theta:S\to\mathbb{R}_{\ge 0}$ is a measurable function satisfying $\int_S\theta\,d\mu\ge 2+\varepsilon$.
	\end{itemize}
	
	We study the following Roth-type inequalities in the variable $\beta\in\mathbb{K}^{\times}$:
	\begin{equation}\label{eq:main_system}
		\left|\beta-\alpha_i\right|_{\omega}\le (AH(\beta))^{-\theta(\omega)},\quad
		\text{for } i=1,\ldots,n \text{ and } \omega\in S_i.
	\end{equation}
	Moreover, when $n=1$, $\Omega$ is equipped with the counting measure, $\theta=2+\varepsilon$, $A=1$, and $S=\{\omega\}$, system \eqref{eq:main_system} reduces to the classical Roth inequality
	$$ |\beta-\alpha|_\omega\le H(\beta)^{-2-\varepsilon}. $$
	
	We first prove a result analogous to \cite[Theorem 1.2]{d}.
	
	\begin{theorem}\label{eq:main_theorem1}
		Fix $A\ge B(\varepsilon)$. Let $\mathcal B$ be the set of solutions $\beta\in\mathbb K^{\times}$ of system \eqref{eq:main_system} satisfying $h(\beta)\ge \log A$. Then
		$$ \#\mathcal B \le (N-1)\left(1+\frac{\log 8nN^2N!}{\log(1+\varepsilon)}\right). $$
	\end{theorem}
	
	Next, we prove the following quantitative result. It gives an upper bound that removes the large-height restriction on the solutions, and the parameter $A$ no longer needs to be bounded below by a large constant.
	
	\begin{theorem}\label{eq:main_theorem2}
		Fix $A\ge 2$. Let $\mathcal C$ be the set of solutions $\beta\in\mathbb K^{\times}$ of system \eqref{eq:main_system}, and put
		$$ T=(1+4/\varepsilon)b(\varepsilon/2)-(4/\varepsilon+2)\log A. $$
		Then
		$$ \#\mathcal C \le 2+(N-1)\left(1+\frac{\log 8nN^2N!}{\log(1+\varepsilon/2)}\right)+\max\left\{1+\frac{\log T-\log c(A)}{\log(1+\varepsilon)},\,0\right\}. $$
	\end{theorem}
	
	Our proof follows the ideas of \cite{dz} and combines three gap principles: the first forward gap principle, the second forward gap principle, and the backward gap principle. The first two are derived from a zeroth forward gap principle analogous to \cite[Proposition 3.1]{d}, but for a different system, while the last is analogous to \cite[Proposition 3.9]{d}, again for a different system. The second forward and first-forward gap principles establish $\varepsilon$-dependent lower bounds for both the difference and the ratio of the heights of solutions. In contrast, the backward gap principle shows that large solutions of system \eqref{eq:main_system} can be contained in several large intervals.
	
	The reason Theorem \ref{eq:main_theorem2} extends the admissible range of $A$ and eliminates the large-height restriction is that we combine the three gap principles with Theorem \ref{eq:main_theorem1}. By applying these four tools iteratively, we obtain a bootstrapping effect on Theorem \ref{eq:main_theorem1}. Our methods partially overlap with those of \cite{bom}; however, the results in Theorem \ref{eq:main_theorem2} are completely new in the setting of adelic curves.
	
	Moreover, in Section \ref{sec:constructive_proof}, we provide a fully constructive quantitative generalization of a special system. Specifically, we consider the systems introduced in \cite[Theorem A]{dz} and \cite[Theorem B]{dz}.
	\begin{theorem}[Theorem A]\label{A}
		Let $\mathbb X=(\mathbb K,\Omega, \phi)$ be a proper adelic curve satisfying the $\mu$-equicontinuity condition and the uniform integrability condition. Fix some distinct elements $\alpha_1,\dots ,\alpha_n\in\mathbb K$. Let $S=\sqcup_{i=1}^{n} S_i \subseteq \Omega$ be a measurable subset of $\Omega$ of finite measure. Then for any $\varepsilon>0$ there exists a real constant $C>0$ (depending on the fixed data) such that for any $\beta\in \mathbb K$  with  $h_{\mathbb X}(\beta)>C$
		it holds that:
		\begin{equation}
			\sum_{i=1}^{n}\int_{S_i}\log^-\vert \beta-\alpha_i\vert_\omega d\mu(\omega)>-(2+\varepsilon)h_\mathbb X(\beta).
		\end{equation}
	\end{theorem}
	\begin{theorem}[Theorem B]\label{B}
		Let $\mathbb X=(\mathbb K,\Omega, \phi)$ be a proper adelic curve satisfying the $\mu$-equicontinuity condition and the uniform integrability condition. Fix some distinct elements $\alpha_1,\dots ,\alpha_n\in\mathbb K$. Let $S$ be a measurable  subset of $\Omega$ of finite measure. Then for any $\varepsilon>0$ and any $c\in\mathbb R$ there exists a real constant $C>0$ (depending on the fixed data) such that for any $\beta\in \mathbb K$  with  $h_{\mathbb X}(\beta)>C$ 
		it holds that:
		\begin{equation}\label{roth systemB}
			\int_{S}\,\min_{1\le i\le n}\left(\log^-\vert \beta-\alpha_i\vert_\omega\right) d\mu(\omega)>-(2+\varepsilon)h_\mathbb X(\beta)+c.
		\end{equation}
	\end{theorem}
	From these two theorems, We analysis the system below
	\begin{equation}\label{roth theoremB1}
		\int_{S}\,\min_{1\le i\le n}\left(\log^-\vert \beta-\alpha_i\vert_\omega\right) d\mu(\omega)\le-(2+\varepsilon)h_\mathbb X(\beta)+c, \quad c \in \mathbb{R}.
	\end{equation}
	The standard Dolce--Vojta method relies on a non-constructive extraction lemma (Lemma \ref{main lemma}) below, which uses the compactness of the median space to extract an infinite convergent subsequence. Consequently, this analytical approach only yields the finiteness of the solution set but fails to provide an explicit, computable bound. 
	
	\begin{lemma}[{\cite[Lemma 8.10]{vo}}]\label{main lemma}
		Let \(\Xi\) be a set, and let \((S,\mathcal{A},\mu)\) be a measure space with
		finite measure. Let \(h:\Xi\to[h_0,\infty)\) be an unbounded function with
		\(h_0>0\). Let \(q\ge 1\) be an integer. For every \(\xi\in\Xi\) and
		\(j=1,\dots,q\), let
		$$
		\lambda_{\xi,j}:S\longrightarrow \mathbb{R}_{\ge 0}
		$$
		be measurable functions. Let \(C_0\ge 0\). Assume:
		
		\begin{enumerate}[label=(H\arabic*),leftmargin=*]
			\item \textbf{(Height integral bound.)} For all \(\xi\in\Xi\) and
			\(j=1,\dots,q\),
			\begin{equation}
				\int_S \lambda_{\xi,j}\,d\mu
				\le h(\xi)+C_0. 
			\end{equation}
			
			\item \textbf{(\(\mu\)-equicontinuity.)} For every \(\delta_{\mathrm{eq}}>0\)
			and every \(\varepsilon_{\mathrm{eq}}>0\) there exists a measurable cover
			\(C_1,\dots,C_m\) of \(S\) such that for each \(\xi\in\Xi\) there is a
			set \(T_\xi\in\mathcal{A}\) with \(\mu(T_\xi)\le\delta_{\mathrm{eq}}\)
			and
			\begin{equation}
				|\lambda_{\xi,j}(\omega)-\lambda_{\xi,j}(\omega')|
				\le \varepsilon_{\mathrm{eq}}\bigl(h(\xi)+C_0\bigr) 
			\end{equation}
			for all \(j=1,\dots,q\), all \(\omega,\omega'\in C_l\setminus T_\xi\),
			and all \(l=1,\dots,m\).
		\end{enumerate}
		
		Then, for every \(\delta_{\mathrm{out}}>0\), there exists a subset
		\(\Xi'\subseteq\Xi\) together with sets \(T_\xi\in\mathcal{A}\)
		\((\xi\in\Xi')\) such that:
		
		\begin{itemize}
			\item \(h\) is unbounded on \(\Xi'\);
			\item \(\mu(T_\xi)\le\delta_{\mathrm{out}}\) for all \(\xi\in\Xi'\);
			\item for all \(\eta,\zeta\in\Xi'\), all
			\(\omega\in S\setminus(T_\eta\cup T_\zeta)\), and all
			\(j=1,\dots,q\),
			\begin{equation}
				\left|
				\frac{\lambda_{\eta,j}(\omega)}{h(\eta)}
				-
				\frac{\lambda_{\zeta,j}(\omega)}{h(\zeta)}
				\right|
				\le
				\left(
				4+\frac{C_0}{h(\eta)}+\frac{C_0}{h(\zeta)}
				\right)\varepsilon_{\mathrm{eq}}.
			\end{equation}
		\end{itemize}
	\end{lemma}
	
	\begin{remark}\label{rem 1.6}
		For adelic curves, the two technical hypotheses of the extraction lemma correspond to natural geometric conditions. A proper adelic curve $\mathbb{X}=(\mathbb{K},\Omega, \mu)$ satisfies the \emph{\(\mu\)-equicontinuity condition} if for any measurable set $S\subset\Omega$ of finite measure and any real numbers $\varepsilon,\delta>0$ there exists a finite measurable cover \(C_1,\ldots, C_m\) of \(S\) satisfying the following: for all $\beta\in\mathbb{K}^\times$ there exists a measurable set \(U_\beta\subset\Omega_\infty\) such that \(\mu(U_\beta)<\delta\) and
		\[
		\left\vert-\log^-\vert\beta\vert_{\omega}+\log^-\vert\beta\vert_{\omega'} \right\vert<\varepsilon h(\beta)\,,\quad \forall \omega,\omega'\in C_j\setminus U_\beta\,,\; \forall j=1,\ldots, m.
		\]
		In \cite[Proposition 8.9]{vo} it is shown that arithmetic function fields satisfy this condition, which provides precisely hypothesis (H2) above. Furthermore, for the control of exceptional sets (often referred to as assumption (H3)), a proper adelic curve $\mathbb{X}=(\mathbb{K},\Omega, \mu)$ satisfies the \emph{uniform integrability condition} if for any \(\varepsilon>0\) there exists \(\delta>0\) such that if $T\subset\Omega_{\infty}$ is a measurable subset satisfying \(\mu(T)<\delta\), then 
		\begin{align}\label{uni_int}
			\int_T -\log^-|\beta|_\omega \, d\mu(\omega)<\varepsilon h_{\mathbb{X}}(\beta)\,,\quad \forall\beta\in\mathbb{K}^\times\,.
		\end{align}

		According to \cite[Proposition 8.8]{vo}, arithmetic function fields also satisfy this condition. 
	\end{remark}
	
	To overcome the limitation of the non-constructive extraction, we replace the continuous compactness with a discrete combinatorial pigeonhole principle. By partitioning the continuous space of median vectors into a finite number of grid boxes, we can bound the number of solutions without invoking limits. The resulting bound will depend explicitly on the topological and measure-theoretic parameters of the system. 
	
	Let \(\mathcal{B}\) denote the set of solutions \(\beta \in \mathbb{K}^\times\) to the system \eqref{roth theoremB1} satisfying \(h(\beta) \ge \log A\). For \(\beta\in\mathcal{B}\) and \(i=1,\dots,n\), define
	\[
	\lambda_{\beta,i}(\omega)
	:=
	-\log^{-}|\beta-\alpha_{i}|_{\omega},
	\qquad
	\omega\in S.
	\]
	By the integral Liouville inequality, hypothesis (H1) of Lemma \ref{main lemma} holds with
	$C_{0}:=\max_{1\le i\le n}h(\alpha_{i})+\log 2.$
	
	We state our constructive result in the following proposition.
	
	\begin{proposition}\label{prop:finite_capacity}
		Let \(0 < \varepsilon' \le \varepsilon\). Let
		\[
		\Gamma=8nN^2N!,\qquad N=N(\frac{\varepsilon'}{2}),
		\]
		and let \(r\) be a real number satisfying
		\[
		\Gamma/2<r\le \Gamma.
		\]
		Assume that \(A> A_1(r,\frac{\varepsilon'}{2})>B(\frac{\varepsilon'}{2})\), where $A_1(r,\frac{\varepsilon'}{2})=\max \left\{{B(\frac{\varepsilon'}{2})^{\left(\frac{4nN^2N!+1}{r-4nN^2N!}\right)}, {B(\frac{\varepsilon'}{2})^{\left(\frac{4}{\varepsilon'}+1\right) }}},e^{h_{\mathrm{min}}} \right\}$, and that hypotheses (H1) and (H2) of Lemma \ref{main lemma}, hypothesis (H3) of Remark \ref{rem 1.6} hold for the functions
		\[
		\lambda_{\beta,i}(\omega)=-\log^-|\beta-\alpha_i|_\omega,
		\qquad \beta\in\mathcal{B},\quad i=1,\dots,n.
		\]
		
		Let \(\delta_{\mathrm{out}}\), \(\varepsilon_{\mathrm{eq}}\) and $h_{\mathrm{mini}}$ (which is a lower height bound defined in (\ref{h_mini})) be parameters depending on \(\varepsilon\), \(\varepsilon'\), \(S\) and the parameters in (H1), (H2) and (H3), which will be specified in Section \ref{sec:constructive_proof}.
		
		Then there exists a finite grid \(\mathcal G\) of boxes and, for every \(\beta\in\mathcal{B}\), an exceptional set \(T_\beta\subset S\) with \(\mu(T_\beta)\le\delta_{\mathrm{out}}\), such that every box in \(\mathcal G\) contains at most \(N-1\) elements of any subset
		$
		\{\beta_1,\dots,\beta_s\}\subseteq\mathcal{B},
		$
		whose consecutive heights satisfy
		$$
		\frac{h(\beta_{j+1})}{h(\beta_j)}\ge r
		\qquad (j=1,\dots,s-1).
		$$
		More precisely, one may take
		$$
		M=\left(\prod_{l=1}^{m}\left\lceil\frac{C_{10,l}}{2\varepsilon_{\mathrm{eq}}}\right\rceil\right)^n
		$$
		boxes, where
		\[
		C_{10,l}:=\frac{2}{\mu(C_l)}\left(1+\frac{C_0}{h_0}\right),
		\qquad h_0:=\log A>0.
		\]
		Consequently, for every subset \(S\subseteq \mathcal{B}\),
		\[
		\#S\le
		\left(1+\frac{\log\Gamma}{\log(1+\varepsilon')}\right)(N-1)M.
		\]
	\end{proposition}
	\section{Preliminaries}
	\label{cha:prel}
	In this section, we review the fundamental definitions of the theory of adelic curves in \cite{Ch-Mo, Ch-Mo2, dz, d}. We use the following notation throughout the article.
	$$
	\log^+ x := \max\{0, \log x\}, \qquad
	\log^- x := \min\{0, \log x\}, \qquad
	\text{for all } x \in \mathbb{R}_{>0}.
	$$
	
	\begin{definition}\label{AdCu}
		Let $\mathbb{K}$ be a field of characteristic $0$, let $M_{\mathbb{K}}$ be the set of all absolute values on $\mathbb{K}$, and let $\Omega=(\Omega,\mathcal{A}, \mu)$ be a measure space endowed with a map
		$$
		\phi \colon \Omega \to M_{\mathbb{K}},\qquad
		\omega \mapsto |\cdot|_\omega := \phi(\omega),
		$$
		such that for any $a\in \mathbb{K}^\times$, the real-valued function $\omega\mapsto \log |a|_\omega$ lies in $L^1(\Omega,\mu)$. The triple $\mathbb{X}=(\mathbb{K},\Omega,\phi)$ is called an \emph{adelic curve}; $\Omega$ and $\phi$ are respectively called the \emph{parameter space} and the \emph{parametrization}. We denote by $\Omega_{\infty}$ the subset of $\Omega$ consisting of all elements $\omega$ for which $|\cdot|_{\omega}$ is Archimedean, and we set $\Omega_0 := \Omega\setminus \Omega_{\infty}$.
	\end{definition}
	
	It is straightforward to verify that the set $\Omega_{\infty}$ is always measurable \cite[Proposition 3.1.1]{Ch-Mo}.
	
	\begin{definition}
		An adelic curve $\mathbb{X}=(\mathbb{K},\Omega,\phi)$ is said to be \emph{proper} if for any $a\in \mathbb{K}^\times$, we have
		\begin{equation}\label{eq:pr_for}
			\int_{\Omega} \log |a|_\omega \, d\mu(\omega) = 0.
		\end{equation}
	\end{definition}
	
	The above property is a natural generalization of the product formula. We briefly mention some important examples of proper adelic curves.
	
	\begin{example}\label{counting}
		Any number field equipped with the counting measure on $\Omega$ and the ``natural'' choice of $\phi$ constitutes a proper adelic curve. Note that in this case, $\phi$ may not be injective; indeed, to obtain a product formula as in equation \eqref{eq:pr_for}, one usually considers the complex embeddings with multiplicity two.
	\end{example}
	
	\begin{example}
		An arithmetic function field $K$ endowed with a big polarization is a proper adelic curve. For further details, see \cite[Section 3]{vo}.
	\end{example}
	
	\begin{definition}
		Let $\mathbb{X}=(\mathbb{K},\Omega,\phi)$ be a proper adelic curve. The \emph{(logarithmic) height} of an element $\beta\in\mathbb{K}^\times$ is defined as
		$$
		h(\beta) = h_{\mathbb{X}}(\beta) = \int_{\Omega}\log^+|\beta|_\omega \, d\mu(\omega).
		$$
		We also define the exponential height as $H(\beta) := e^{h(\beta)}$.
	\end{definition}
	
	For $\omega\in\Omega_\infty$, Ostrowski's theorem guarantees the existence of a real number $\epsilon(\omega)\in (0,1]$ such that $\vert\cdot\vert_\omega = \vert\cdot\vert^{\epsilon(\omega)}$, where the right-hand side denotes the standard Euclidean absolute value on $\mathbb{C}$. Thus we obtain a map $\epsilon\colon\Omega_{\infty}\to (0,1]$, which can be extended to $\epsilon\colon\Omega\to [0,1]$ by setting $\epsilon|_{\Omega_0} := 0$. For instance, for an Archimedean absolute value $\vert\cdot\vert_\omega$, we have $\log \vert 2\vert_\omega = \epsilon(\omega)\log 2$, which yields the explicit expression
	$$
	\epsilon(\omega) = \frac{\log^+\vert 2 \vert_\omega}{\log 2}
	$$
	on the entire space $\Omega_\infty$. Clearly, $\epsilon$ is a measurable function. Moreover, we usually normalize the measure on $\Omega$ so that $h(2)\le \log 2$.
	
	If instead $h(2)=d\log 2$ with $d\ge1$, then, after replacing $|\cdot|_\omega$ by $|\cdot|_\omega^{1/d}$ and passing to the normalized height $H^{1/d}$, the system becomes
	$$
	\left|\beta-\alpha_i\right|^{1/d}_{\omega} \le (A^{1/d}H(\beta)^{1/d})^{-\theta(\omega)},
	\qquad \text{for } i=1,\ldots, n \text{ and } \omega\in S_i.
	$$
	Since $H(2)^{1/d}=2$, we can therefore work with the normalized adelic curve and the new parameter $A^{1/d}$, obtaining a new system that satisfies $h(2)\le \log 2$.
	
	One can prove several standard properties for the height on $\mathbb{X}$, but for our purposes we will only require the following.
	
	\begin{proposition}[Integral Liouville inequality]\label{prop:liou}
		For any $\alpha,\beta\in\mathbb{K}^\times$ and any measurable subset $U\subseteq \Omega$, we have
		$$
		\int_U\log \vert \alpha-\beta\vert_\omega \, d\mu(\omega)
		\ge \int_U\log^- \vert \alpha-\beta\vert_\omega \, d\mu(\omega)
		\ge -\log 2 - h(\alpha) - h(\beta).
		$$
	\end{proposition}
	
	\begin{proof}
		See \cite[Proposition 1.13]{dz}.
	\end{proof}
	
	\section{Gap principles}\label{gap}
	
	\subsection{The forward gap principles}
	In this section, we establish $\varepsilon$-dependent lower bounds for both the difference and the ratio of the heights of solutions. We first state a proposition whose proof is analogous to that of \cite[Proposition 3.1]{d}.
	
	Recall the definition of $c(X)$ in Section \ref{cha_intro}.
	
	\begin{proposition}\label{amb_gap}[Zeroth forward gap principle]
		Let $\beta_1,\beta_2\in\mathbb K^\times$ be two distinct solutions of the system of inequalities \eqref{eq:main_system} such that
		\begin{equation}\label{eq:cond_h_amb}
			h(\beta_2)\ge h(\beta_1).
		\end{equation}
		Then
		$$
		h(\beta_2)\ge (1+\varepsilon)h(\beta_1)+c(A).
		$$
	\end{proposition}
	
	\begin{proof}
		Set $S_0=S\cap \Omega_0$ and $S_\infty=S\cap \Omega_\infty$. For $\omega\in S$, there exists $i\in\{1,\dots,n\}$ such that $\omega\in S_i$.
		
		First suppose that $\omega$ is Archimedean. Recall the definition of $\epsilon(\omega)$ in Section \ref{cha:prel}. By the triangle inequality,
		$$
		|\beta_1-\beta_2|^{1/\epsilon(\omega)}_\omega
		\le |\beta_1-\alpha_i|^{1/\epsilon(\omega)}_\omega
		+|\beta_2-\alpha_i|^{1/\epsilon(\omega)}_\omega
		\le 2(AH(\beta_1))^{-\theta(\omega)/\epsilon(\omega)},
		$$
		which is independent of the choice of $i$. Hence
		$$
		|\beta_1-\beta_2|_\omega\le 2^{\epsilon(\omega)}(AH(\beta_1))^{-\theta(\omega)}.
		$$
		Taking logarithms and integrating over $S_\infty$, we obtain
		\begin{equation}\label{eps1}
			\int_{S_\infty}-\log|\beta_1-\beta_2|_\omega
			\ge (\log A+h(\beta_1))\int_{S_\infty}\theta\,d\mu
			-\log 2\int_{S_\infty}\epsilon\,d\mu
			\ge (\log A+h(\beta_1))\int_{S_\infty}\theta\,d\mu-h(2).
		\end{equation}
		
		Next suppose that $\omega$ is non-Archimedean. Then
		$$
		|\beta_1-\beta_2|_\omega
		\le \max\{|\beta_1-\alpha_i|_\omega,|\beta_2-\alpha_i|_\omega\}
		\le (AH(\beta_1))^{-\theta(\omega)}.
		$$
		After integrating over $S_0$, this gives
		\begin{equation}\label{eps2}
			\int_{S_0}-\log|\beta_1-\beta_2|_\omega
			\ge (\log A+h(\beta_1))\int_{S_0}\theta\,d\mu.
		\end{equation}
		Adding \eqref{eps1} and \eqref{eps2} and using $\int_S\theta\,d\mu\ge 2+\varepsilon$, we get
		$$
		\int_S-\log|\beta_1-\beta_2|_\omega
		\ge (\log A+h(\beta_1))\int_S\theta\,d\mu-h(2)
		\ge (2+\varepsilon)(\log A+h(\beta_1))-h(2).
		$$
		By Proposition \ref{prop:liou},
		$$
		h(\beta_1)+h(\beta_2)+\log 2
		\ge (2+\varepsilon)(\log A+h(\beta_1))-h(2).
		$$
		Rearranging the inequality yields
		$$
		h(\beta_2)\ge (1+\varepsilon)h(\beta_1)+(2+\varepsilon)\log A-h(2)-\log 2.
		$$
		Since $h(2)\le \log 2$, the conclusion follows.
	\end{proof}
	
	While the above gap principle estimates the separation of heights, concrete bounds require expressions that specify the direction of estimation. We therefore specialize the gap principle above to deduce the following two immediate and key corollaries.
	
	\begin{corollary}\label{wave_gap}[First forward gap principle]
		Assume that $A\ge2$. Let $\beta_1,\beta_2\in\mathbb K^\times$ be two distinct solutions of the system of inequalities \eqref{eq:main_system} such that
		\begin{equation}\label{eq:cond_h_wave}
			h(\beta_2)\ge h(\beta_1)>0.
		\end{equation}
		Then
		$$
		\frac{h(\beta_2)}{h(\beta_1)}\ge1+\varepsilon.
		$$
	\end{corollary}
	
	\begin{proof}
		By Proposition \ref{amb_gap} and $A\ge2$, we have $c(A)\ge0$. Hence
		$$
		h(\beta_2)\ge (1+\varepsilon)h(\beta_1).
		$$
		where the last inequality uses $h(\beta_1)>0$.
	\end{proof}
	
	\begin{corollary}\label{part_gap}[Second forward gap principle]
		Assume that $A\ge2$. Then the system \eqref{eq:main_system} has at most one solution in the interval $[r,r+c(A))$ for every $r\ge0$.
	\end{corollary}
	
	\begin{proof}
		Suppose that $\beta_1,\beta_2\in\mathbb K^\times$ are two distinct solutions of \eqref{eq:main_system} such that
		\begin{equation}\label{eq:cond_h_part}
			h(\beta_2)\ge h(\beta_1)\ge r.
		\end{equation}
		Then Proposition \ref{amb_gap} gives
		\begin{align*}
			h(\beta_2)
			&\ge (1+\varepsilon)h(\beta_1)+c(A)\\
			&\ge (1+\varepsilon)r+c(A)\\
			&\ge r+c(A),
		\end{align*}
		a contradiction.
	\end{proof}

	\subsection{The backward gap principle}
	\begin{remark}
		The numerical constraints on $A$ and $h(\beta)$ in this section differ from those in the previous section.
	\end{remark}
	
	We begin with some notions from \cite{dz} and \cite{d}. The key theorem below is taken from \cite{dz}.
	
	\begin{definition}\label{d3.5}
		Let $I$ be a bounded real interval with $a=\inf I$, $b=\sup I$ and $1<a<b$. The \emph{logarithmic length} of $I$ is defined as the real number $\log\frac{b}{a}$. If $a<1$, one can still define the logarithmic length after an adequate translation of the interval.
	\end{definition}
	
	\begin{definition}
		We say that two vectors $\bm{\alpha}=(\alpha^{(1)},\ldots,\alpha^{(N)})$ and $\bm{\beta}=(\beta^{(1)},\ldots,\beta^{(N)})$ are \emph{componentwise different} if $\alpha^{(j)}\neq\beta^{(j)}$ for $j=1,\ldots,N$.
	\end{definition}
	
	Consider a proper adelic curve $(\mathbb K,\Omega,\phi)$ and vectors $\bm{\alpha}_1,\ldots,\bm{\alpha}_n,\bm{\beta}\in\mathbb K^N$ which are componentwise different. We construct the following matrices $T:=T(\bm{\alpha}_1,\ldots,\bm{\alpha}_n)\in M(n\times N,\mathbb K)$ and $T(\bm{\beta})\in M((n+1)\times N,\mathbb K)$:
	$$
	T:=\begin{pmatrix}
		\alpha^{(1)}_1 & \alpha^{(2)}_1 & \ldots & \alpha^{(N)}_1 \\
		\alpha^{(1)}_2 & \alpha^{(2)}_2 & \ldots & \alpha^{(N)}_2 \\
		\vdots & \vdots & \vdots & \vdots \\
		\alpha^{(1)}_n & \alpha^{(2)}_n & \ldots & \alpha^{(N)}_n
	\end{pmatrix},
	\qquad
	T(\bm{\beta})=\begin{pmatrix}
		\alpha^{(1)}_1 & \alpha^{(2)}_1 & \ldots & \alpha^{(N)}_1 \\
		\alpha^{(1)}_2 & \alpha^{(2)}_2 & \ldots & \alpha^{(N)}_2 \\
		\vdots & \vdots & \vdots & \vdots \\
		\alpha^{(1)}_n & \alpha^{(2)}_n & \ldots & \alpha^{(N)}_n \\
		\beta^{(1)} & \beta^{(2)} & \ldots & \beta^{(N)}
	\end{pmatrix}.
	$$
	We denote by $\bm{\alpha}^{(j)}$, $j=1,\ldots,N$, the columns of $T$, that is,
	$$
	\bm{\alpha}^{(j)}=\begin{pmatrix}
		\alpha^{(j)}_1 \\
		\alpha^{(j)}_2 \\
		\vdots \\
		\alpha^{(j)}_n
	\end{pmatrix},
	$$
	and by $\bm{\alpha}_{h}$, for $h\in\{1,2,\dots,n\}$, the rows:
	$$
	\bm{\alpha}_{h}=\begin{pmatrix}
		\alpha^{(1)}_h, \alpha^{(2)}_h,\ldots, \alpha^{(N)}_h
	\end{pmatrix}.
	$$
	
	\begin{definition}\label{boxing}
		For the matrix $T(\bm{\beta})$, consider the following real numbers for any $j=1,\ldots,N$:
		\begin{equation}\label{pivot}
			\rho_{j}:=4^{2N!}H(\beta^{(j)})\prod_{h=1}^{n}H(\alpha_h^{(j)})^{\frac{2N!}{n}},
		\end{equation}
		\begin{equation}\label{pivotdue}
			\rho'_{j}:=4^{N!}H(\beta^{(j)})\prod_{h=1}^{n}H(\alpha_h^{(j)})^{\frac{N!}{2n}}.
		\end{equation}
		We say that $T(\bm{\beta})$ satisfies the \emph{$h$-gap property} if
		$$
		\frac{\log \rho_j}{\log \rho'_{j+1}}<\frac{1}{4nN^2N!},
		\qquad \forall j=1,\ldots,N-1.
		$$
	\end{definition}
	
	\begin{definition}\label{bounding}
		Fix a measurable set $S=S_1\sqcup\ldots\sqcup S_n$. We say that an integrable function $\theta\colon S\to\mathbb R_{\ge 0}$ is a \emph{column bounding function for} $T(\bm{\beta})$ on $S$ if
		\begin{equation}\label{boundingbis}
			-\frac{1}{\log \rho_{j}}\log\vert \alpha^{(j)}_h-\beta^{(j)}\vert_\omega\ge\theta(\omega)
			\qquad
			\forall j=1,\ldots,N,\ \forall h=1,\ldots,n,\ \forall\omega\in S_h.
		\end{equation}
	\end{definition}
	
	\begin{remark}\label{rem_forthelemma}
		If $\theta$ is a column bounding function, then it immediately follows that
		$$
		\vert \alpha^{(j)}_h-\beta^{(j)}\vert_{\omega}\le 1
		$$
		for any $j=1,\dots,N$, $h=1,\dots,n$, and $\omega\in S_h$.
	\end{remark}
	
	The following nontrivial theorem is the crucial ingredient for the second gap principle.
	
	\begin{theorem}\label{thm:teoremadue}
		Let $(\mathbb K,\Omega,\phi)$ be a proper adelic curve. Assume that the matrix $T(\bm{\beta})\in M((n+1)\times N,\mathbb K)$ satisfies the $h$-gap property (see Definition \ref{boxing}). Moreover, fix $N>21^2\log 2n$ and let $\theta\colon S\to\mathbb R_{\ge 0}$ be a column bounding function for $T(\bm{\beta})$ (see Definition \ref{bounding}). Then
		$$
		\int_{S}\theta\,d\mu<2+\frac{7\sqrt{\log 2n}}{\sqrt{N}}.
		$$
	\end{theorem}
	
	\begin{proof}
		See \cite[Theorem 3.9]{dz}.
	\end{proof}
	
	The following gap principle, which is an analog of \cite[Proposition 3.9]{d}, informally says that for $A$ large enough, one cannot allow too many large solutions of \eqref{eq:main_system} whose heights are very spaced.
	
	\begin{proposition}\label{prop:large_gap}[Backward gap principle]
		Fix $A\ge B(\varepsilon)$. Then the solutions $\beta\in\mathbb K^\times$ of \eqref{eq:main_system} such that $H(\beta)\ge A$ are contained in at most $N-1$ intervals of logarithmic length less than or equal to $\Gamma=8nN^2N!$.
	\end{proposition}
	
	\begin{proof}
		Assume the contrary. Then the set of solutions with \(H(\beta)\ge A\) cannot be covered by \(N-1\) intervals of logarithmic length $\Gamma$. By the usual greedy covering argument on the logarithmic height axis, there exists
		$$
		\beta_1,\dots,\beta_N\in\mathcal{C},
		\qquad h(\beta_1)<\dots<h(\beta_N),
		$$
		such that
		\begin{equation*}
			h(\beta_{j+1})\notin[h(\beta_j),\Gamma h(\beta_j)]
		\end{equation*}
		So  
		\begin{equation}\label{eq:separated}
			\frac{h(\beta_{j+1})}{h(\beta_j)}> \Gamma
			\qquad (j=1,\dots,N-1).
		\end{equation}
		
		We claim that \(T(\beta_1,\dots,\beta_N)\) satisfies the \(h\)-gap property of Definition~\ref{boxing}. Indeed, from the definitions of $\rho_j$ and $\rho'_{j+1}$ in \eqref{pivot} and \eqref{pivotdue}, there is a constant \(C_1\), depending only on \(n,N,\alpha_1,\dots,\alpha_n\), such that
		$$
		\log\rho_j\le h(\beta_j)+C_1,
		\qquad
		\log\rho'_{j+1}\ge h(\beta_{j+1})-C_1 .
		$$
		Using \eqref{eq:separated} and \(A\ge B(\varepsilon)\), we may absorb \(C_1\) and obtain
		$$
		\frac{\log\rho_j}{\log\rho'_{j+1}}
		\le
		\frac{h(\beta_j)+C_1}{{\Gamma}h(\beta_j)}
		<
		\frac{1}{4nN^2N!},
		$$
		for all $j=1,\dots,N-1$. Hence, the \(h\)-gap property holds.
		
		Next, we show that \(\theta\) is a column-bounding function for \(T(\beta_1,\dots,\beta_N)\). Let \(\omega\in S_h\). Since each \(\beta_j\) solves \eqref{eq:main_system}, we have
		$$
		-\log|\beta_j-\alpha_h|_\omega
		\ge
		\theta(\omega)\bigl(\log A+h(\beta_j)\bigr).
		$$
		Thus
		$$
		-\frac{1}{h(\beta_j)+\log A}\log|\beta_j-\alpha_h|_\omega
		\ge
		\theta(\omega).
		$$
		Since \(A\ge B(\varepsilon)\) absorbs the constants in $\rho_j$, we have
		$$
		\log\rho_j\le h(\beta_j)+\log A,
		$$
		and therefore
		$$
		-\frac{1}{\log\rho_j}\log|\beta_j-\alpha_h|_\omega
		\ge
		\theta(\omega).
		$$
		So \(\theta\) is column-bounding.
		
		By Theorem~\ref{thm:teoremadue},
		$$
		\int_S\theta\,d\mu
		<
		2+\frac{7\sqrt{\log(2n)}}{\sqrt N}.
		$$
		Since \(N=N(\varepsilon)\), we have
		$$
		\frac{7\sqrt{\log(2n)}}{\sqrt N}\le\varepsilon.
		$$
		Consequently,
		$$
		2+\varepsilon
		\le
		\int_S\theta\,d\mu
		<
		2+\varepsilon,
		$$
		a contradiction. Therefore, solutions with \(H(\beta)\ge A\) are contained in at most \(N-1\) intervals of logarithmic length at most $\Gamma$.
	\end{proof}
	
	\section{Proof of Theorem \ref{eq:main_theorem2}}
	
	We first prove an important lemma.
	
	\begin{lemma}\label{log_interval}
		Fix $A\ge 2$. Let $\mathcal C$ be the set of solutions $\beta\in\mathbb K^\times$ of system \eqref{eq:main_system}. If $I$ is an interval of logarithmic length less than or equal to $\Gamma$, then
		$$
		\#\bigl(H(\mathcal{C})\cap I\bigr)\le \frac{\log \Gamma}{\log (1+\varepsilon)}+1.
		$$
	\end{lemma}
	
	\begin{proof}
		Let
		$$
		\inf I\le x_1\le x_2\le\ldots\le x_r\le \sup I
		$$
		be a finite sequence of real numbers such that $x_j^{1+\varepsilon}\le x_{j+1}$. Clearly, we have
		$$
		(1+\varepsilon)^{r-1}\le\frac{\log x_r}{\log x_1}\le\Gamma,
		$$
		which means that
		$$
		r\le 1+\frac{\log\Gamma}{\log(1+\varepsilon)}.
		$$
		Using Corollary \ref{wave_gap}, we obtain
		$$
		\#\bigl(H(\mathcal{C})\cap I\bigr)\le \frac{\log \Gamma}{\log (1+\varepsilon)}+1.
		$$
	\end{proof}
	
	\textbf{Proof of Theorem \ref{eq:main_theorem1}.}
	Since $A\ge B(\varepsilon)\ge4^2>2$, we have $c(A)>0$. Combined with Corollary \ref{part_gap}, this shows that the height function $h$ is injective on $\mathcal B$. Therefore
	$$
	\#\mathcal{B}=\#h(\mathcal{B})=\#H(\mathcal{B}).
	$$
	Hence, by Proposition \ref{prop:large_gap} and Lemma \ref{log_interval},
	$$
	\#\mathcal{B}=\#H(\mathcal{B})
	\le (N-1)\left(1+\frac{\log 8nN^2N!}{\log(1+\varepsilon)}\right).
	$$
	\hfill{$\square$}
	
	\begin{remark}
		Fix $A\ge 2$. Let $\mathcal C$ be the set of solutions $\beta\in\mathbb K^\times$ of system \eqref{eq:main_system}. Then $h$ is also injective on $\mathcal{C}$. Indeed,
		$$
		c(A)\ge c(2)=(2+\varepsilon)\log2-2\log 2=\varepsilon\log2>0,
		$$
		and the injectivity follows from Corollary \ref{part_gap}.
	\end{remark}
	
	Before proving Theorem \ref{eq:main_theorem2}, we need the following lemma.
	
	\begin{lemma}\label{eq:main_theorem1.5}
		Fix $A\ge B(\varepsilon)$. Let $\mathcal C_l$ be the set of solutions $\beta\in\mathbb K^\times$ of system \eqref{eq:main_system}. Then
		$$
		\#\mathcal C_l \le (N-1)\left(1+\frac{\log 8nN^2N!}{\log(1+\varepsilon)}\right)+1.
		$$
	\end{lemma}
	
	\begin{proof}
		Since
		$$
		c(A)-\log A=(1+\varepsilon)\log A-2\log2>0,
		$$
		Corollary \ref{part_gap} gives
		$$
		\#\bigl(h(\mathcal{C}_l)\cap[0,\log A)\bigr)
		\le \#\bigl(h(\mathcal{C}_l)\cap[0,c(A))\bigr)\le 1.
		$$
		By Theorem \ref{eq:main_theorem1}, we obtain
		\begin{align*}
			\#\mathcal{C}_l
			&=\#h(\mathcal{C}_l)\\
			&= \#\bigl(h(\mathcal{C}_l)\cap[0,\log A)\bigr)
			+\#\bigl(h(\mathcal{C}_l)\cap[\log A,\infty)\bigr)\\
			&\le (N-1)\left(1+\frac{\log 8nN^2N!}{\log(1+\varepsilon)}\right)+1.
		\end{align*}
	\end{proof}
	
	\textbf{Proof of Theorem \ref{eq:main_theorem2}.}
	For $\beta\in \mathbb{K}^{\times}$ satisfying
	
	$$
	h(\beta)\ge (1+4/\varepsilon)b(\varepsilon/2)-(4/\varepsilon+2)\log A,
	$$
	we have
	\begin{align*}
		h(\beta)&\ge (1+4/\varepsilon)b(\varepsilon/2)-(4/\varepsilon+2)\log A\\
		\implies H(\beta)&\ge\frac{B(\varepsilon/2)^{1+4/\varepsilon}}{A^{4/\varepsilon+2}}\\
		\implies (B(\varepsilon/2)H(\beta))^{-\frac{2+\varepsilon/2}{2+\varepsilon}\theta(\omega)}
		&\ge (AH(\beta))^{\theta(\omega)}.
	\end{align*}
	Since
	$$
	\int_S\frac{2+\varepsilon/2}{2+\varepsilon}\theta(\omega)\,d\mu\ge 2+\varepsilon/2,
	$$
	Lemma \ref{eq:main_theorem1.5} yields
	$$
	\#\bigl(h(\mathcal{C})\cap[(1+4/\varepsilon)b(\varepsilon/2)-(4/\varepsilon+2)\log A,\infty)\bigr)
	\le (N-1)\left(1+\frac{\log 8nN^2N!}{\log(1+\varepsilon/2)}\right)+1.
	$$
	
	We distinguish two cases.
	
	\underline{Case-$\mathbb{I}$} If
	$$
	\log\bigl[(1+4/\varepsilon)b(\varepsilon/2)-(4/\varepsilon+2)\log A\bigr]-\log c(A)\le0,
	$$
	then
	\begin{align*}
		\#\mathcal{C}
		&=\#h(\mathcal{C})\\
		&\le \#\bigl(h(\mathcal{C})\cap[(1+4/\varepsilon)b(\varepsilon/2)-(4/\varepsilon+2)\log A,\infty)\bigr)\\
		&\quad+\#\bigl(h(\mathcal{C})\cap[0,c(A))\bigr)\\
		&\le (N-1)\left(1+\frac{\log 8nN^2N!}{\log(1+\varepsilon/2)}\right)+2\\
		&\le 2+(N-1)\left(1+\frac{\log 8nN^2N!}{\log(1+\varepsilon/2)}\right)\\
		&\quad+\max\left\{1+\frac{\log\bigl[(1+4/\varepsilon)b(\varepsilon/2)-(4/\varepsilon+2)\log A\bigr]-\log c(A)}{\log(1+\varepsilon)},\,0\right\}.
	\end{align*}
	
	\underline{Case-$\mathbb{II}$} If
	$$
	\log\bigl[(1+4/\varepsilon)b(\varepsilon/2)-(4/\varepsilon+2)\log A\bigr]-\log c(A)>0,
	$$
	then, together with Lemma \ref{log_interval}, we deduce
	\begin{align*}
		\#\mathcal{C}
		&=\#h(\mathcal{C})\\
		&\le \#\bigl(h(\mathcal{C})\cap[(1+4/\varepsilon)b(\varepsilon/2)-(4/\varepsilon+2)\log A,\infty)\bigr)\\
		&\quad+\#\bigl(h(\mathcal{C})\cap[0,c(A))\bigr)\\
		&\quad+\#\left(H(\mathcal{C})\cap\left[e^{c(A)},\frac{B(\varepsilon/2)^{1+4/\varepsilon}}{A^{4/\varepsilon+2}}\right)\right)\\
		&\le (N-1)\left(1+\frac{\log 8nN^2N!}{\log(1+\varepsilon/2)}\right)+2\\
		&\quad+\frac{\log\left(\log\left(\frac{B(\varepsilon/2)^{1+4/\varepsilon}}{A^{4/\varepsilon+2}}\right)/\log(e^{c(A)})\right)}{\log(1+\varepsilon)}\\
		&\le 2+(N-1)\left(1+\frac{\log 8nN^2N!}{\log(1+\varepsilon/2)}\right)\\
		&\quad+\max\left\{1+\frac{\log\bigl[(1+4/\varepsilon)b(\varepsilon/2)-(4/\varepsilon+2)\log A\bigr]-\log c(A)}{\log(1+\varepsilon)},\,0\right\}.
	\end{align*}
	In both cases, let $$ T=(1+4/\varepsilon)b(\varepsilon/2)-(4/\varepsilon+2)\log A. $$
	This completes the proof.

	\section{A constructive quantitative generalization of system (\ref{roth theoremB1})}
	\label{sec:constructive_proof}
	
	In the proof, we mainly rely on Vojta’s arithmetic function field version (See \cite{vo}). According to the discussion in \cite[Theorem B]{dz}, its technical approach is the same as that used to prove the adelic curve version, and it can 
	be directly generalized to the setting of adelic curve. In the main-line proof (which relies on arithmetic function fields), we will intersperse some lemmas and definitions for adelic curves, and we should know that this conversion is legitimate. We first recall the system (\ref{roth theoremB1}) and a lemma introduced in Section \ref{cha_intro}.
	
	\begin{theorem}[Theorem \ref{A}]
		Let $\mathbb X=(\mathbb K,\Omega, \phi)$ be a proper adelic curve satisfying the $\mu$-equicontinuity condition and the uniform integrability condition. Fix some distinct elements $\alpha_1,\dots ,\alpha_n\in\mathbb K$. Let $S=\sqcup_{i=1}^{n} S_i \subseteq \Omega$ be a measurable subset of $\Omega$ of finite measure. Then for any $\varepsilon>0$ there exists a real constant $C>0$ (depending on the fixed data) such that for any $\beta\in \mathbb K$  with  $h_{\mathbb X}(\beta)>C$ 
		it holds that:
		\begin{equation*}
			\sum_{i=1}^{n}\int_{S_i}\log^-\vert \beta-\alpha_i\vert_\omega d\mu(\omega)>-(2+\varepsilon)h_\mathbb X(\beta).
		\end{equation*}
	\end{theorem}
	\begin{theorem}[Theorem \ref{B}]
		Let $\mathbb X=(\mathbb K,\Omega, \phi)$ be a proper adelic curve satisfying the $\mu$-equicontinuity condition and the uniform integrability condition. Fix some distinct elements $\alpha_1,\dots ,\alpha_n\in\mathbb K$. Let $S$ be a measurable  subset of $\Omega$ of finite measure. Then for any $\varepsilon>0$ and any $c\in\mathbb R$ there exists a real constant $C>0$ (depending on the fixed data) such that for any $\beta\in \mathbb K$  with  $h_{\mathbb X}(\beta)>C$ 
		it holds that:
		\begin{equation*}
			\int_{S}\,\min_{1\le i\le n}\left(\log^-\vert \beta-\alpha_i\vert_\omega\right) d\mu(\omega)>-(2+\varepsilon)h_\mathbb X(\beta)+c.
		\end{equation*}
	\end{theorem}

	According to these theorems, we have system (\ref{roth theoremB1})
	\begin{equation*}
		\int_{S_i}\min_{1\le i\le n}\left(\log^-\vert \beta-\alpha_i\vert_\omega\right) d\mu(\omega)\le-(2+\varepsilon)h_\mathbb X(\beta)+c, \quad c \in \mathbb{R}.
	\end{equation*}
	
	\begin{lemma}[Lemma \ref{main lemma}]\label{lem:extraction}
		Let \(\Xi\) be a set, and let \((S,\mathcal{A},\mu)\) be a measure space with
		finite measure. Let \(h:\Xi\to[h_0,\infty)\) be an unbounded function with
		\(h_0>0\). Let \(q\ge 1\) be an integer. For every \(\xi\in\Xi\) and
		\(j=1,\dots,q\), let
		$$
		\lambda_{\xi,j}:S\longrightarrow \mathbb{R}_{\ge 0}
		$$
		be measurable functions. Let \(C_0\ge 0\). Assume:
		
		\begin{enumerate}[label=(H\arabic*),leftmargin=*]
			\item \textbf{(Height integral bound.)} For all \(\xi\in\Xi\) and
			\(j=1,\dots,q\),
			\begin{equation}\label{H1}
				\int_S \lambda_{\xi,j}\,d\mu
				\le h(\xi)+C_0. 
			\end{equation}
			
			\item \textbf{(\(\mu\)-equicontinuity.)} For every \(\delta_{\mathrm{eq}}>0\)
			and every \(\varepsilon_{\mathrm{eq}}>0\) there exists a measurable cover
			\(C_1,\dots,C_m\) of \(S\) such that for each \(\xi\in\Xi\) there is a
			set \(T_\xi\in\mathcal{A}\) with \(\mu(T_\xi)\le\delta_{\mathrm{eq}}\)
			and
			\begin{equation}\label{H2}
				|\lambda_{\xi,j}(\omega)-\lambda_{\xi,j}(\omega')|
				\le \varepsilon_{\mathrm{eq}}\bigl(h(\xi)+C_0\bigr) 
			\end{equation}
			for all \(j=1,\dots,q\), all \(\omega,\omega'\in C_l\setminus T_\xi\),
			and all \(l=1,\dots,m\).
		\end{enumerate}
		
		Then, for every \(\delta_{\mathrm{out}}>0\), there exists a subset
		\(\Xi'\subseteq\Xi\) together with sets \(T_\xi\in\mathcal{A}\)
		\((\xi\in\Xi')\) such that:
		
		\begin{itemize}
			\item \(h\) is unbounded on \(\Xi'\);
			\item \(\mu(T_\xi)\le\delta_{\mathrm{out}}\) for all \(\xi\in\Xi'\);
			\item for all \(\eta,\zeta\in\Xi'\), all
			\(\omega\in S\setminus(T_\eta\cup T_\zeta)\), and all
			\(j=1,\dots,q\),
			\begin{equation} \label{H3}
				\left|
				\frac{\lambda_{\eta,j}(\omega)}{h(\eta)}
				-
				\frac{\lambda_{\zeta,j}(\omega)}{h(\zeta)}
				\right|
				\le
				\left(
				4+\frac{C_0}{h(\eta)}+\frac{C_0}{h(\zeta)}
				\right)\varepsilon_{\mathrm{eq}}.
			\end{equation}
		\end{itemize}
	\end{lemma}
	For \(\beta\in\mathcal B\) and \(i=1,\dots,n\), define
	\[
	\lambda_{\beta,i}(\omega)
	:=
	-\log^{-}|\beta-\alpha_{i}|_{\omega},
	\qquad
	\omega\in S.
	\]
	By the integral Liouville inequality, hypothesis (H1) of
	Lemma \ref{lem:extraction} holds with
	\[
	C_{0}:=\max_{1\le i\le n}h(\alpha_{i})+\log 2.
	\]
	Now we give our main result.
	\begin{proposition}\label{main 2}
		Let \(0 < \varepsilon' \le \varepsilon\). Let
		\[
		\Gamma=8nN^2N!,\qquad N=N(\frac{\varepsilon'}{2}),
		\]
		and let \(r\) be a real number satisfying
		\[
		\Gamma/2<r\le \Gamma.
		\]
		Assume that \(A > A_1(r,\frac{\varepsilon'}{2})>B(\frac{\varepsilon'}{2})\), where $A_1(r,\frac{\varepsilon'}{2})=\max \left\{{B(\frac{\varepsilon'}{2})^{\left(\frac{4nN^2N!+1}{r-4nN^2N!}\right)}, {B(\frac{\varepsilon'}{2})^{\left(\frac{4}{\varepsilon'}+1\right) }}},e^{h_{\mathrm{min}}}\right\}$, and that hypotheses (H1) and (H2) of Lemma \ref{lem:extraction}, hypothesis (H3) of Remark \ref{rem 1.6} hold for the functions
		\[
		\lambda_{\beta,i}(\omega)=-\log^-|\beta-\alpha_i|_\omega,
		\qquad \beta\in\mathcal{B},\quad i=1,\dots,n.
		\]
		
		Let \(\delta_{\mathrm{out}}\), \(\varepsilon_{\mathrm{eq}}\) and $h_{\mathrm{min}}$ (see equation (\ref{h_mini})) be parameters depending on \(\varepsilon\) , \(\varepsilon'\) , \(S\) and the parameters in (H1), (H2) and (H3).
		
		Then there exists a finite grid \(\mathcal G\) of boxes and, for every \(\beta\in\mathcal{B}\), an exceptional set \(T_\beta\subset S\) with \(\mu(T_\beta)\le\delta_{\mathrm{out}}\), such that every box in \(\mathcal G\) contains at most \(N-1\) elements of any subset
		\[
		\{\beta_1,\dots,\beta_s\}\subseteq\mathcal{B}
		\]
		whose consecutive heights satisfy
		\[
		\frac{h(\beta_{j+1})}{h(\beta_j)}\ge r
		\qquad (j=1,\dots,s-1).
		\]
		More precisely, one may take
		\[
		M=\left(\prod_{l=1}^{m}\left\lceil\frac{C_{10,l}}{2\varepsilon_{\mathrm{eq}}}\right\rceil\right)^n
		\]
		boxes, where
		\[
		C_{10,l}:=\frac{2}{\mu(C_l)}\left(1+\frac{C_0}{h_0}\right),
		\qquad h_0:=\log A>0.
		\]
		Consequently, for every subset \(S\subseteq \mathcal{B}\),
		\[
		\#S\le
		\left(1+\frac{\log\Gamma}{\log(1+\varepsilon')}\right)(N-1)M.
		\]
	\end{proposition}
	
	\begin{proof}
		
		We first construct the grid. Let \(\delta_{\mathrm{eq}}=\delta_{\mathrm{out}}/3n\). According to hypothesis (H2), there is a finite measurable cover \(C_1,\dots,C_m\) of \(S\) such that for every \(\beta\in\mathcal B\) there exists sets \(T_{\beta,j}^0\subset S\) with
		\[
		\mu(T_{\beta,j}^0)\le\delta_{\mathrm{eq}}
		\]
		and
		\[
		|\lambda_{\beta,j}(\omega)-\lambda_{\beta,j}(\omega')|
		\le \varepsilon_{\mathrm{eq}}\bigl(h(\beta)+C_0\bigr)
		\]
		for all \(\omega,\omega'\in C_l\setminus T_{\beta,j}^0\), and all \(l=1,\dots,m\)

		For each \(\beta\in\mathcal B\) and each \(l\) with \(\mu(C_l)>0\), define the median
		\[
		m_{\beta,j,l}
		:=
		\inf\left\{
		t\in\mathbb R:
		\mu\bigl(\{\omega\in C_l:\lambda_{\beta,j}(\omega)\ge t h(\beta)\}\bigr)
		\le \frac{\mu(C_l)}2
		\right\},
		\]
		
		Put
		\[
		V_{\beta,j,l}:=\{\omega\in C_l:\lambda_{\beta,j}(\omega)\le m_{\beta,j,l}h(\beta)\},
		\qquad
		W_{\beta,j,l}:=\{\omega\in C_l:\lambda_{\beta,l}(\omega)\ge m_{\beta,l}h(\beta)\}.
		\]
		If $\mu(C_l \cap T_{\beta,j}^0) < \mu(C_l)/2$, using (H1) and the definition of \(W_{\beta,j,l}\), we obtain
		\[
		\frac{\mu(C_l)}2 m_{\beta,j,l}h(\beta)
		\le
		\int_{W_{\beta,j,l}\setminus T_\beta^0}\lambda_{\beta,j}\,d\mu
		\le
		\int_S\lambda_{\beta,j}\,d\mu
		\le
		h(\beta)+C_0.
		\]
		Hence
		\[
		m_{\beta,j,l}
		\le
		C_{10,l}
		:=
		\frac{2}{\mu(C_l)}\left(1+\frac{C_0}{h_0}\right).
		\]
		For each \(\beta\), define a vector \(m_{\beta,j}\in\mathbb R^m\) by setting its \(l\)-th coordinate equal to \(m_{\beta,j,l}\) if $\mu(C_l \cap T_{\beta,l}^0) < \mu(C_l)/2$, and equal to \(0\) otherwise. Then
		\[
		m_{\beta,j}\in\prod_{l=1}^m[0,C_{10,l}].
		\] 
		and
		\[
		(m_{\beta,j})_{j=1}^n\in\prod_{j=1}^n\prod_{l=1}^m[0,C_{10,l}]
		\]
		Partition each interval \([0,C_{10,l}]\) into subintervals of length \(2\varepsilon_{\mathrm{eq}}\). This gives a grid with
		\[
		M=\prod_{j=1}^n\prod_{l=1}^m\left\lceil\frac{C_{10,l}}{2\varepsilon_{\mathrm{eq}}}\right\rceil
		\]
		boxes. 
		
		Assign each \(\beta\) to the box containing $(m_{\beta,j})_{j=1}^n$. Define $T_{\beta,j} = T_{\beta,j}^0 \cup \bigcup_{\mu(C_l \cap T_{\beta,j}^0) \ge \mu(C_l)/2} C_l$. Then $\mu(T_{\beta,j}) \le 3\mu(T_{\beta,j}^0) \le \delta_{\mathrm{out}}/n$. Let  $$T_\beta=\bigcup_{i=1}^n T_{\beta,j}.$$ We have \(\mu(T_\beta)\le\delta_{\mathrm{out}}\).
		
		Now, we compare two elements \(\beta,\gamma\) in the same box. For \(\omega\in C_l\setminus(T_{\beta}\cup T_{\gamma})\), by the construction of $T_\beta$ and, we have $\mu(C_l \cap T_\beta^0) < \mu(C_l)/2$. Thus, there exist $\omega_1 \in C_l \setminus T_\beta$ such that
		\[
		\lambda_{\beta,j}(\omega)
		\ge
		\lambda_{\beta,j}(\omega_1)-\epsilon_{\mathrm{eq}}\bigl(h(\beta)+C_0\bigr)
		\ge
		m_{\beta,j,l}h(\beta)-\epsilon_{\mathrm{eq}}\bigl(h(\beta)+C_0\bigr).
		\]
		Choosing \(\omega_0\in V_{\beta,j,l}\setminus T_\beta^0\), we similarly get
		\[
		\lambda_{\beta,j}(\omega)
		\le
		m_{\beta,l}h(\beta)+\epsilon_{\mathrm{eq}}\bigl(h(\beta)+C_0\bigr).
		\]
		Therefore
		\[
		\left|
		\frac{\lambda_{\beta,j}(\omega)}{h(\beta)}-m_{\beta,j,l}
		\right|
		\le
		\varepsilon_{\mathrm{eq}}\left(1+\frac{C_0}{h(\beta)}\right).
		\]
		The same estimate holds for \(\gamma\). Since \(\beta\) and \(\gamma\) lie in the same box,
		\[
		|m_{\beta,j,l}-m_{\gamma,j,l}|\le 2\varepsilon_{\mathrm{eq}}.
		\]
		Combining these estimates yields
		\[
		\left|
		\frac{\lambda_{\beta,j}(\omega)}{h(\beta)}
		-\frac{\lambda_{\gamma,j}(\omega)}{h(\gamma)}
		\right|
		\le
		\left(4+\frac{C_0}{h(\beta)}+\frac{C_0}{h(\gamma)}\right)\varepsilon_{\mathrm{eq}}.
		\]
		Because \(h(\beta),h(\gamma)\ge h_{\mathrm{min}}\), we obtain the uniform bound
		\begin{align}\label{box_principle}
			\left|
			\frac{\lambda_{\beta,i}(\omega)}{h(\beta)}
			-\frac{\lambda_{\gamma,i}(\omega)}{h(\gamma)}
			\right|
			\le
			\left(4+\frac{2C_0}{h_{\mathrm{min}}}\right)\varepsilon_{\mathrm{eq}}
			=:\varepsilon_{12}  
		\end{align}

		for all \(i=1,\dots,n\) and all \(\omega\in S\setminus(T_\beta\cup T_\gamma)\).

		We also require the following lemmas to extract the appropriate sequences aimed at deriving a contradiction from the local box bound. 
		
		We choose $c' \in \mathbb{R}$ such that $c' \ge \log 4$ throughout this paper.
		
		We choose $h_{\min} \in \mathbb{R}$ such that
		\begin{align}\label{h_mini}
			h_{\min} > {\frac{c'+\max\{-c,0\}}{\varepsilon-\varepsilon'}}
		\end{align} throughout this paper.

		\begin{lemma} \label{prop:Gap}
			[Quantitive version of \cite[Proposition 8.12]{vo}]
			Let $\mathbb K$ be a polarized arithmetic function field. Assume that 
			\begin{equation}\label{roth_encore}
				\int_{S}\,\min_{1\le i\le n}\left(\log^-\vert \beta-\alpha_i\vert_\omega\right) d\mu(\omega)>-(2+\varepsilon)h_\mathbb X(\beta)+c,\quad c \in \mathbb{R}.
			\end{equation}
			
			Let $r_0\ge h_{\mathrm{min}}>0 ,\quad r_1\ge1$ and
			\[
			B:=\left\{
			\beta\in\mathbb{K}:
			\int_{S}\min_{1\le i\le n}
			\bigl(\log^{-}|\beta-\alpha_i|_{\omega}\bigr)\,d\mu
			\le -(2+\varepsilon_0)h(\beta)+c
			\right\}.
			\]

			Define
			\[
			t_0:=\left\lceil\frac{\log r_1}{\log(1+\varepsilon)}\right\rceil.
			\]
			If
			\[
			\#\{\beta\in B:h(\beta)\ge r_0\}\ge t_0(N-1)M+1,
			\]
			for certain $S$, $\alpha_1,\ldots,\alpha_n\in\mathbb K$, $\varepsilon_0\in\mathbb R_{>0}$ and $c \in \mathbb{R}$.
			
			Then there exists $\beta_1,\dots,\beta_{N}\in B$ such that:
			\begin{enumerate}
				\item[(1)] $h(\beta_1)\ge r_0\ge h_{\mathrm{min}}>0$;
				\item[(2)] 
				\[
				\frac{h(\beta_k)}{h(\beta_{k-1})}\ge r_1\ge1
				\qquad (k=2,\dots,N);
				\]
				\item[(3)] there exists a partition
				$S=S_1\sqcup\cdots\sqcup S_n$ such that
				\begin{equation}\label{5.1}
					\sum_{i=1}^{n}\int_{S_i}
					\min_{1\le k\le N}
					\left(
					\frac{-\log^{-}|\beta_k-\alpha_i|_{\omega}}{h(\beta_k)}
					\right)d\mu
					\ge2+\varepsilon'+\frac{c'}{h_\X(\xi_1)},
				\end{equation}
				
				where $0<\varepsilon'<\varepsilon$ is a given small quantity.
			\end{enumerate}
		\end{lemma}
		\begin{proof}
			Let \[\Xi=\{\beta\in B:h(\beta)\ge r_0\}\]
			
			Define
			$$
			\lambda_{\xi,j}(\omega)
			:=
			-\log^-|\xi-\alpha_j|_\omega,
			\qquad
			h:=h_\X.
			$$
			By the integral Liouville inequality, hypothesis (H1) of
			Lemma~\ref{lem:extraction} holds with
			$
			C_0=\max_j h_\X(\alpha_j)+\log 2.
			$
			
			The \(\mu\)-equicontinuity condition (established in  \cite[Proposition 8.9]{vo} ) provides hypothesis (H2), and the
			uniform integrability condition provides hypothesis (H3) (established in \cite[Lemma 8.8]{vo} ) below.
			
			We use the following consequence, which is the
			gap-sequence extraction step. Under (H1), (H2), and the uniform integrability
			condition
			\begin{enumerate}[label=(H3),leftmargin=*]
				\item for every \(\varepsilon_{\mathrm{int}}>0\) there exists
				\(\delta_{\mathrm{int}}>0\) such that for every measurable set
				\(T\subseteq S\) with \(\mu(T)\le\delta_{\mathrm{int}}\) and every
				\(\xi\in\Xi\), \(j=1,\dots,q\),
				$$
				\int_T\lambda_{\xi,j}\,d\mu
				\le
				\varepsilon_{\mathrm{int}}\bigl(h(\xi)+C_0\bigr),
				$$
			\end{enumerate}
			for every \(N\ge 1\), \(\varepsilon_{\mathrm{gap}}>0\), \(c_{\mathrm{gap}}>0\),
			there exist \(\xi_1,\dots,\xi_N\in\Xi\) (Corollary \ref{amb_gap}), a subset
			\(T\subseteq S\), and a measurable function
			$$
			J:S\setminus T\longrightarrow\{1,\dots,n\}
			$$
			such that
			\begin{equation}
				h(\xi_1)\ge r_0\ge h'_{\min}>0
			\end{equation}
			\begin{equation}\label{5.4}
				\frac{h(\xi_i)}{h(\xi_{i-1})}\ge r_1
				\qquad (i=2,\dots,N), 
			\end{equation}
			and
			\begin{equation}\label{5.5}
				\int_{S\setminus T}
				\min_{1\le i'\le N}
				\frac{\lambda_{\xi_{i'},J(\omega)}(\omega)}{h(\xi_{i'})}
				\,d\mu(\omega)
				\ge
				\int_S\frac{\lambda_{\xi_i}}{h(\xi_i)}\,d\mu
				-\varepsilon_{\mathrm{gap}}
				+\frac{c_{\mathrm{gap}}}{h(\xi_1)}, 
			\end{equation}
			where \(\lambda_{\xi_i}=\max_{1\le j\le n}\lambda_{\xi_i,j}\) and $h'_{\min}\in\mathbb{R}$ is a chosen constant satisfied 
			$$h'_{\min}>\frac{c_{\mathrm{gap}}}{\varepsilon_{\mathrm{gap}}}$$.
			
			We now prove this extraction step. Choose $\varepsilon_{\mathrm{int}}>0$, $\varepsilon_{\mathrm{eq}}>0$, and $ h'_{\mathrm{min}}>0$ such that
			\begin{equation}\label{h_min}
				n\varepsilon_{\mathrm{int}}
				\left(1+\frac{C_0}{h'_{\min}}\right)
				+
				\left(4+\frac{2C_0}{h'_{\min}}\right)\varepsilon_{\mathrm{eq}}\mu(S)
				+
				\frac{c_{\mathrm{gap}}}{h'_{\min}}
				\le
				\varepsilon_{\mathrm{gap}}.
			\end{equation}
			
			Using the parameters $\delta_{\mathrm{out}}=\delta_{\mathrm{int}}/N$ and $\varepsilon_{\mathrm{eq}}$, we can construct a grid consisting of $M$ boxes. We obtain the corresponding $\Xi'$ in one selected box, a subset
			\(\Xi'_1\subseteq\Xi\) and sets \(T_\xi\) with
			\(\mu(T_\xi)\le\delta_{\mathrm{out}}\), such that (\ref{box_principle}) holds. Using the pigeonhole principle, we can let
			\(\#\Xi'_1\ge t_0(N-1)+1\).
			Write its elements in increasing order of height:
			\[
			h(\beta_1)<h(\beta_2)<\cdots<h(\beta_s),
			\qquad s:=\#\Xi'_1.
			\]
			Since \(\beta_1,\dots,\beta_s\) are elements of \(\Xi'_1\), another version of first-forward gap principle is still satisfied for \(\Xi'\).
			It gives, for every \(j=1,\dots,s-1\),
			\[
			h(\beta_{j+1})\ge(1+\varepsilon')h(\beta_j).
			\]
			
			Choose
			\(\xi_1,\dots,\xi_N\in\Xi'_1\) with
			$$
			\frac{h(\xi_i)}{h(\xi_{i-1})}\ge r_1
			\qquad (i=2,\dots,N),
			$$
			and with \(\xi_1\) of minimal height among them. 
			Indeed, $\xi_i=\beta_{(i-1)t_0+1}$ will be satisfied.
			
			Let
			$$
			T=\bigcup_{i=1}^N T_{\xi_i}.
			$$
			Then
			$$
			\mu(T)\le N\delta_{\mathrm{out}}\le\delta_{\mathrm{int}}.
			$$
			Because \(\lambda_{\xi_i}\le\sum_{j=1}^q\lambda_{\xi_i,j}\), hypothesis (H3)
			yields
			\begin{equation}\label{5.6}
				\int_T\lambda_{\xi_i}\,d\mu
				\le
				\sum_{j=1}^n\int_T\lambda_{\xi_i,j}\,d\mu
				\le
				n\varepsilon_{\mathrm{int}}\bigl(h(\xi_i)+C_0\bigr)
				\le
				n\varepsilon_{\mathrm{int}}
				\left(1+\frac{C_0}{h'_{\min}}\right)h(\xi_i). 
			\end{equation}
			
			For \(\omega\in S\setminus T\), define the normalized functions
			$$
			\lambda_{ij}=\frac{\lambda_{\xi_i,j}(\omega)}{h(\xi_i)}.
			$$
			Let
			$$
			\varepsilon_{12}
			=
			\left(4+\frac{2C_0}{h(\xi_1)}\right)\varepsilon_{\mathrm{eq}}.
			$$
			By (\ref{H3}), we know
			$$
			|\lambda_{ij}-\lambda_{i'j}|\le\varepsilon_{12}
			\qquad
			\text{for all } i,i',j.
			$$
			We claim that there is some index \(j\in\{1,\dots,q\}\) such that
			\begin{equation}\label{5.7}
				\max_{1\le j'\le q}\lambda_{ij'}
				-
				\min_{1\le i'\le N}\lambda_{i'j}
				\le
				\varepsilon_{12}
				\qquad
				\text{for all } i=1,\dots,N.
			\end{equation}
			To see this, let \((i^*,j^*)\) be the pair that maximizes
			\(\lambda_{i^*j^*}\). Pick \(j=j^*\). By the bound on differences,
			$$
			\max_{j'}\lambda_{ij'}
			\le
			\lambda_{i^*j^*}
			\le
			\min_{i'}\lambda_{i'j^*}+\varepsilon_{12},
			$$
			proving (\ref{5.7}). Let \(J(\omega)\) be the smallest such \(j\), ensuring that \(J\) is
			measurable.
			
			Now, substituting these bounds:
			$$
			\begin{aligned}
				&\int_S\frac{\lambda_{\xi_i}}{h(\xi_i)}\,d\mu
				-
				\int_{S\setminus T}
				\min_{1\le i'\le N}
				\frac{\lambda_{\xi_{i'},J(\omega)}(\omega)}{h(\xi_{i'})}
				\,d\mu
				+
				\frac{c_{\mathrm{gap}}}{h(\xi_1)}\\
				&=
				\int_T\frac{\lambda_{\xi_i}}{h(\xi_i)}\,d\mu
				+
				\int_{S\setminus T}
				\left(
				\max_{1\le j'\le n}\lambda_{ij'}
				-
				\min_{1\le i'\le N}\lambda_{i',J(\omega)}
				\right)d\mu
				+
				\frac{c_{\mathrm{gap}}}{h(\xi_1)}\\
				&\le
				q\varepsilon_{\mathrm{int}}
				\left(1+\frac{C_0}{h_{\min}}\right)
				+
				\varepsilon_{12}\mu(S\setminus T)
				+
				\frac{c_{\mathrm{gap}}}{h'_{\min}}
				\le
				\varepsilon_{\mathrm{gap}}.
			\end{aligned}
			$$
			This proves (\ref{5.5}).
			
			Now apply the extraction step with
			$$
			\varepsilon_{\mathrm{gap}}=\varepsilon-\varepsilon',
			\qquad
			c_{\mathrm{gap}}=c'+\max\{-c,0\},
			\qquad
			h'_{\min}=h_{\min}
			$$
			and gap ratio \(r_1\) and a lower bound $r_0$. We obtain \(\xi_1,\dots,\xi_N\in\Xi\), a subset
			\(T\subseteq S\), and a function
			$$
			J:S\setminus T\to\{1,\dots,n\}
			$$
			satisfying (\ref{5.4}) and (\ref{5.5}).
			
			Let \(\lambda_\xi=\max_{1\le j\le n}\lambda_{\xi,j}\). Substituting (\ref{5.5})
			into (\ref{5.1}) gives
			$$
			\begin{aligned}
				&\int_{S\setminus T}
				\min_{1\le i'\le N}
				\frac{\lambda_{\xi_{i'},J(\omega)}(\omega)}{h_\X(\xi_{i'})}
				\,d\mu(\omega)\\
				&\ge
				\int_S\frac{\lambda_{\xi_i}}{h_\X(\xi_i)}\,d\mu
				-\varepsilon_{\mathrm{gap}}
				+\frac{c_{\mathrm{gap}}}{h_\X(\xi_1)}\\
				&>
				\left(
				2+\varepsilon+\frac{c}{h_\X(\xi_i)}
				\right)
				-
				(\varepsilon-\varepsilon')
				+
				\frac{c'+\max\{-c,0\}}{h_\X(\xi_1)}\\
				&\ge
				2+\varepsilon'+\frac{c'}{h_\X(\xi_1)}.
			\end{aligned}
			$$
			Let \(S_j=J^{-1}(j)\) for \(j=1,\dots,n\). Since \(J\) is defined on
			\(S\setminus T\), the subsets \(S_j\) are disjoint, and summing the integrals
			over \(S_j\) yields (\ref{5.5}).
		\end{proof}
		\begin{remark} For adelic curves, the rather technical proof is identical to Lemma \ref{prop:Gap} for arithmetic function fields. 
			We explain why the assumptions hold immediately in adelic curves: assumption (H1) follows by the integral Liouville inequality if we put \(\Xi=\{\beta_k\}\), \(q=n\), \(\lambda_{\beta_k,j}=-\log^-\vert\alpha_j-\beta_k\vert\) and \(C_0= \log 2+\max_j h(\alpha_j)\) for \(j=1,\ldots, n\). Assumption (H2) is the \(\mu\)-equicontinuity condition. Assumption (H3) is the uniform integrability condition. 
		\end{remark}
		Thus, we state the following quantitative lemma over adelic curves, which is consistent with Lemma \ref{prop:Gap}.
		\begin{lemma}[Quantitative finite extraction over adelic curves]\label{vo8.12}
			Let $\mathbb X=(\mathbb K,\Omega, \phi)$ be a proper adelic curve satisfying the $\mu$-equicontinuity condition  and the uniform integrability condition. Assume that 
			\begin{equation}\label{roth_encore}
				\int_{S}\,\min_{1\le i\le n}\left(\log^-\vert \beta-\alpha_i\vert_\omega\right) d\mu(\omega)>-(2+\varepsilon)h_\mathbb X(\beta)+c.
			\end{equation}
			
			Let $r_0\ge h_{\mathrm{min}}>0,\quad r_1\ge1$ and
			\[
			B:=\left\{
			\beta\in\mathbb{K}:
			\int_{S}\min_{1\le i\le n}
			\bigl(\log^{-}|\beta-\alpha_i|_{\omega}\bigr)\,d\mu
			\le -(2+\varepsilon)h(\beta)+c
			\right\}.
			\]

			Define
			\[
			t_0:=\left\lceil\frac{\log r_1}{\log(1+\varepsilon')}\right\rceil.
			\]
			If
			\[
			\#\{\beta\in B:h(\beta)\ge r_0\}\ge t_0(N-1)M+1,
			\]
			for certain $S$, $\alpha_1,\ldots,\alpha_n\in\mathbb K$, $\varepsilon\in\mathbb R_{>0}$ and $c\in \mathbb{R}$.
			
			Then there exists $\beta_1,\dots,\beta_{N}\in B$ such that:
			\begin{enumerate}
				\item[(1)] $h(\beta_1)\ge r_0\ge h_{\mathrm{min}}>0$;
				\item[(2)] 
				\[
				\frac{h(\beta_k)}{h(\beta_{k-1})}\ge r_1\ge 1
				\qquad (k=2,\dots,N);
				\]
				\item[(3)] there exists a partition
				$S=S_1\sqcup\cdots\sqcup S_n$ such that
				\[
				\sum_{i=1}^{n}\int_{S_i}
				\min_{1\le k\le N}
				\left(
				\frac{-\log^{-}|\beta_k-\alpha_i|_{\omega}}{h(\beta_k)}
				\right)d\mu
				\ge2+\varepsilon'+\frac{c'}{h_\X(\xi_1)},
				\]
				where $0<\varepsilon'<\varepsilon$ is a given small quantity.
			\end{enumerate}
		\end{lemma}
		
		We next prove the result of Proposition \ref{main 2}.
		
		For the case \(n=1\), it follows immediately from the integral Liouville inequality that any approximant \(\beta\) satisfying
		\[
		h(\beta)>\frac{\log 2+h(\alpha_1)+c}{1+\varepsilon}
		\]
		also satisfies the desired inequality.
		
		Now fix \(n\ge 2\), and assume by contradiction that
		\[
		\#S>
		\left(1+\frac{\log\Gamma}{\log(1+\varepsilon')}\right)(N-1)M.
		\]
		
		We apply Lemma \ref{vo8.12}; thus, for $r_0= \log A\ge h_{\min}$, $\Gamma/2< r_1=r\le \Gamma $ and \(0<\varepsilon'<\varepsilon\), there exist \(\beta_1,\ldots,\beta_N\in\mathcal{B}\) satisfying conditions (1)-(3) of the lemma.
		Consider the matrix constructed from the vector \(\bm{\beta}=(\beta_1,\ldots,\beta_N)\):
		\[
		T(\bm{\beta})=\begin{pmatrix}
			\alpha_1& \alpha_1  &\ldots &\alpha_1 \\
			\alpha_2& \alpha_2  &\ldots &\alpha_2\\
			\vdots & \vdots  & \vdots & \vdots\\
			\alpha_n& \alpha_n  &\ldots &\alpha_n\\
			\beta_1& \beta_2  &\ldots &\beta_N
		\end{pmatrix}.
		\] 
		
		For \(j=1,\dots,N\), define
		\[
		\rho_j
		:=
		4^{2N!}H(\beta_j)\prod_{h=1}^n H(\alpha_h)^{\frac{2N!}{n}}.
		\]
		
		Choose 
		\[
		L=\log A^{\frac{\varepsilon'}{4+\varepsilon'}}>\log B(\frac{\varepsilon'}{2})= \log\left(4^{2N!}\prod^n_{i=1}H(\alpha_i)^{\frac{2N!}{n}}\right) 
		\]       
		which will serve as a threshold. 
		
		We now construct a function \(\theta\in L^1(S,\mu)\) that will yield the desired contradiction. Define \(\theta\) as a piecewise function by setting, for each \(\omega\in S_i\),
		\[
		\theta|_{S_i}(\omega)=\min_k\frac{-\log^- \vert \beta_k-\alpha_i\vert_{\omega}}{L+h(\beta_k)}= \min_k \left(\frac{-\log^- \vert \beta_k-\alpha_i\vert_{\omega}}{h(\beta_k)}A_k\right),
		\]
		where $$A_k=1-\frac{L}{L+h(\beta_k)}.$$
		Next, define the set
		\[
		\hat S_i:=\{\omega\in S_i\colon \vert\beta_k-\alpha_i\vert_\omega < 1 \text{ for } k=1,\ldots, N\},
		\]
		and put \(\hat S=\bigcup^n_{i=1} \hat S_i\). Note that by definition, \(\theta|_{S\setminus\hat S}=0\).  Consequently, we can ensure that
		\begin{equation}\label{theepss}
			\varepsilon'>\frac{2}{A_1}+\frac{7\sqrt{\log 2n}}{A_1\sqrt N}-2>0, \qquad \frac{\varepsilon'}{2} \ge \frac{7\sqrt{\log 2n}}{\sqrt N}.
		\end{equation}
		Also, we have $$A_1=1-\frac{L}{L+h(\beta_1)}\ge 1-\frac{\frac{\varepsilon'}{2}}{2+\varepsilon'}=1-\frac{1}{\frac{4}{\varepsilon'}+2}$$ 
		
		Using property (3) of Lemma \ref{vo8.12} and \eqref{theepss}, we obtain:
		\begin{eqnarray}\label{final_contr1}
			\int_{\hat S}\theta \, d\mu =\int_{ S} \theta \, d\mu=\sum^n_{i=1}\int_{S_i}\theta|_{S_i} \, d\mu &\ge& A_1\sum_{i=1}^n\int_{S_i}\min_{k}\frac{-\log^-|\beta_k-\alpha_i|_\omega}{h(\beta_k)}d\mu(\omega) \nonumber \\
			&\ge& A_1\left(2+\varepsilon'\right)\ge 2+\frac{\varepsilon'}{2}\ge2+\frac{7\sqrt{\log 2n}}{\sqrt N}\,.
		\end{eqnarray}
		
		Finally, we verify the conditions for applying Theorem \ref{thm:teoremadue}:
		\begin{itemize}
			\item[\(\bullet\)] The matrix \(T(\bm{\beta})\) satisfies the \(h\)-gap condition:
			
			By the definition of \(B(\frac{\varepsilon'}{2})\) and \(A\ge A_0(r,\frac{\varepsilon'}{2})=e^{\left(C_1\frac{4nN^2N!+1}{r-4nN^2N!}\right)} >B(\frac{\varepsilon'}{2})\), there is a constant \(C_1=\log B(\frac{\varepsilon'}{2})\), depending only on \(n,N,\alpha_1,\dots,\alpha_n\), such that
			\[
			\log\rho_j\le h(\beta_j)+C_1,
			\qquad
			\log\rho_{j+1}\ge h(\beta_{j+1}).
			\]
			Set \(\delta:=r-\Gamma/2>0\). 
			
			Since \(r>\Gamma/2=4nN^2N!\) and \(h(\beta_j)\ge\log A\) with \(A> A_0(r,\frac{\varepsilon'}{2})=e^{\left(C_1\frac{4nN^2N!+1}{r-4nN^2N!}\right)}\) sufficiently large, so that \(\log A> C_1(\frac{4nN^2N!+1}{\delta})\), we may absorb \(C_1\) and obtain
			\[
			\frac{\log\rho_j}{\log\rho_{j+1}}
			\le
			\frac{h(\beta_j)+C_1}{r h(\beta_j)}
			<
			\frac{1}{4nN^2N!}
			\]
			for all \(j=1,\dots,N-1\). Hence the matrix satisfies the \(h\)-gap property of Definition \ref{boxing}.
			\item[\(\bullet\)] Since \(\vert\alpha_i-\beta_k\vert_\omega < 1\) for \(\omega\in \hat S_i\), the function \(\theta\) is column bounding for \(T(\bm{\beta})\) on \(\hat S\). Indeed, for any \(k=1,\ldots, N\) and \(\omega\in \hat S_i\), we have
			\[
			\theta(\omega)\le \frac{-\log^- \vert \beta_k-\alpha_i\vert_{\omega}}{L+h(\beta_k)}\le \frac{-\log^{-}\vert\beta_k-\alpha_i\vert_\omega}{\log\left(4^{2N!}\prod^n_{i=1}H(\alpha_i)^{\frac{2N!}{n}}\right) + h(\beta_k)}.
			\]
		\end{itemize}
		Therefore, all hypotheses of Theorem \ref{thm:teoremadue} are met, allowing us to conclude that
		$$\int_{\hat S}\theta \, d \mu<2+\frac{7\sqrt{\log 2n}}{\sqrt{N}},$$
		which directly contradicts inequality \eqref{final_contr1}, completes the proof.
		
		Therefore, 
		\[
		\#S\le
		\left(1+\frac{\log\Gamma}{\log(1+\varepsilon')}\right)(N-1)M.
		\]
		
		Throughout this paper, we proceed by contradiction with respect to the given quantities. For the detailed counting part, see our Appendix \ref{counting part}.
	\end{proof}
	
	\section{Conclusion}
	
	Throughout this paper, the most significant and technically demanding result is Theorem~\ref{eq:main_theorem2} and Proposition \ref{main 2}. Nevertheless, we believe that our contribution extends beyond these theorems, since the underlying technique has considerable potential for generalization.

	More precisely, the constructive part in Proposition \ref{main 2} replaces nonconstructive compactness extraction with a discrete pigeonhole argument and a finite grid, yielding an explicit finite-capacity bound. This bound depends on measure-theoretic covering, equicontinuity, and uniform integrability parameters, and may be numerically large. Future work may optimize these parameters and extend the framework to more general adelic curve systems and higher-dimensional quantitative approximation problems.

	\appendix
	\section{First forward gap principle for system ($\ref{roth theoremB1}$) (in $\Xi'_1$)}\label{AppendixA}
	
	Recall that we choose $c' \in \mathbb{R}$ such that $c' \ge \log 4$. 
	
	First we give an theorem similar to the Lemma \ref{vo8.12}
	
	\begin{lemma}[Box partition lemma]\label{lem:box_partition}
		Let \(1\le K\le N\), for any \(K\) elements \(\xi_1,\dots,\xi_K\in \Xi'_1\) ordered so that
		\[
		h(\xi_1)\le h(\xi_2)\le\cdots\le h(\xi_K),
		\]
		put \(T:=\bigcup_{i=1}^K T_{\xi_i}\). Then \(\mu(T)\le K\delta_{\mathrm{out}}\le\delta_{\mathrm{int}}\), and there exists a measurable function
		\[
		J:S\setminus T\longrightarrow\{1,\dots,n\}
		\]
		such that for every \(i=1,\dots,K\),
		\begin{equation}\label{5.5}
			\int_{S\setminus T}
			\min_{1\le i'\le K}
			\frac{\lambda_{\xi_{i'},J(\omega)}(\omega)}{h(\xi_{i'})}
			\,d\mu(\omega)
			\ge
			\int_S\frac{\lambda_{\xi_i}}{h(\xi_i)}\,d\mu
			-\varepsilon_{\mathrm{gap}}
			+\frac{c_{\mathrm{gap}}}{h(\xi_1)}
			\ge
			2+\varepsilon'+\frac{c'}{h(\xi_1)},
		\end{equation}
		where \(\lambda_{\xi_i}:=\max_{1\le j\le n}\lambda_{\xi_i,j}=-\min_{1\le j\le n}\log^-|\xi_i-\alpha_j|_\omega\). 
		
		In particular, setting \(S_j:=J^{-1}(j)\) for \(j=1,\dots,n\) defines a measurable partition \(S\setminus T=S_1\sqcup\cdots\sqcup S_n\) satisfying
		\begin{equation}\label{5.1}
			\sum_{j=1}^{n}\int_{S_j}
			\min_{1\le k\le K}
			\left(
			\frac{-\log^{-}|\xi_k-\alpha_j|_{\omega}}{h(\xi_k)}
			\right)d\mu(\omega)
			\ge2+\varepsilon'+\frac{c'}{h(\xi_1)}.
		\end{equation}
	\end{lemma}
	\begin{proof}
		Since \(K\le N\) and \(\mu(T_{\xi_i})\le\delta_{\mathrm{out}}=\delta_{\mathrm{int}}/N\), we have
		\[
		\mu(T)\le K\delta_{\mathrm{out}}\le N\delta_{\mathrm{out}}=\delta_{\mathrm{int}}.
		\]
		Because $$\lambda_{\xi_i}\le\sum_{j=1}^n\lambda_{\xi_i,j}$$ and \(h(\xi_i)\ge\log A\ge h_{\min}\), hypothesis (H3) yields
		\begin{equation}\label{5.6}
			\int_T\lambda_{\xi_i}\,d\mu
			\le
			\sum_{j=1}^n\int_T\lambda_{\xi_i,j}\,d\mu
			\le
			n\varepsilon_{\mathrm{int}}\bigl(h(\xi_i)+C_0\bigr)
			\le
			n\varepsilon_{\mathrm{int}}
			\left(1+\frac{C_0}{h_{\min}}\right)h(\xi_i). 
		\end{equation}
		For \(\omega\in S\setminus T\), define the normalized functions
		\[
		\lambda_{ij}:=\frac{\lambda_{\xi_i,j}(\omega)}{h(\xi_i)}
		\qquad (1\le i\le K,\ 1\le j\le n),
		\]
		and let
		\[
		\varepsilon_{12}
		:=
		\left(4+\frac{2C_0}{h(\xi_1)}\right)\varepsilon_{\mathrm{eq}}
		\le
		\left(4+\frac{2C_0}{h_{\min}}\right)\varepsilon_{\mathrm{eq}}.
		\]
		Since \(\xi_1,\dots,\xi_K\in \Xi'_1\) in the same box, the box proximity estimate ensures that
		\[
		|\lambda_{ij}-\lambda_{i'j}|\le\varepsilon_{12}
		\qquad
		\text{for all } i,i'\in\{1,\dots,K\}\text{ and } j\in\{1,\dots,n\}.
		\]
		We claim that for each \(\omega\in S\setminus T\) there is some index \(j\in\{1,\dots,n\}\) such that
		\begin{equation}\label{5.7}
			\max_{1\le j'\le n}\lambda_{ij'}
			-
			\min_{1\le i'\le K}\lambda_{i'j}
			\le
			\varepsilon_{12}
			\qquad
			\text{for all } i=1,\dots,K.
		\end{equation}
		To see this, let \((i^*,j^*)\in\{1,\dots,K\}\times\{1,\dots,n\}\) be a pair that maximizes \(\lambda_{i^*j^*}\), and pick \(j=j^*\). By the bound on differences,
		\[
		\max_{1\le j'\le n}\lambda_{ij'}
		\le
		\lambda_{i^*j^*}
		\le
		\min_{1\le i'\le K}\lambda_{i'j^*}+\varepsilon_{12},
		\]
		proving \eqref{5.7}. Let \(J(\omega)\in\{1,\dots,n\}\) be the smallest such \(j\), ensuring that \(J:S\setminus T\to\{1,\dots,n\}\) is measurable.
		
		Combining \eqref{5.6}, \eqref{5.7}, and \eqref{h_min}, we obtain
		\[
		\begin{aligned}
			&\int_S\frac{\lambda_{\xi_i}}{h(\xi_i)}\,d\mu
			-
			\int_{S\setminus T}
			\min_{1\le i'\le K}
			\frac{\lambda_{\xi_{i'},J(\omega)}(\omega)}{h(\xi_{i'})}
			\,d\mu
			+
			\frac{c_{\mathrm{gap}}}{h(\xi_1)}\\
			&=
			\int_T\frac{\lambda_{\xi_i}}{h(\xi_i)}\,d\mu
			+
			\int_{S\setminus T}
			\left(
			\max_{1\le j'\le n}\lambda_{ij'}
			-
			\min_{1\le i'\le K}\lambda_{i',J(\omega)}
			\right)d\mu
			+
			\frac{c_{\mathrm{gap}}}{h(\xi_1)}\\
			&\le
			n\varepsilon_{\mathrm{int}}
			\left(1+\frac{C_0}{h_{\min}}\right)
			+
			\varepsilon_{12}\mu(S\setminus T)
			+
			\frac{c_{\mathrm{gap}}}{h_{\min}}
			\le
			\varepsilon_{\mathrm{gap}}.
		\end{aligned}
		\]
		Since \(\xi_i\in\mathcal{B}\) satisfies system \eqref{roth theoremB1}, we have \(\int_S \lambda_{\xi_i}\,d\mu \ge (2+\varepsilon)h(\xi_i)-c\). Recalling that \(\varepsilon_{\mathrm{gap}}=\varepsilon-\varepsilon'\) and \(c_{\mathrm{gap}}=c'+\max\{-c,0\}\), we deduce
		\[
		\begin{aligned}
			\int_{S\setminus T}
			\min_{1\le i'\le K}
			\frac{\lambda_{\xi_{i'},J(\omega)}(\omega)}{h(\xi_{i'})}
			\,d\mu(\omega)
			&\ge
			\int_S\frac{\lambda_{\xi_i}}{h(\xi_i)}\,d\mu
			-\varepsilon_{\mathrm{gap}}
			+\frac{c_{\mathrm{gap}}}{h(\xi_1)}\\
			&\ge
			\left(
			2+\varepsilon-\frac{c}{h(\xi_i)}
			\right)
			-
			(\varepsilon-\varepsilon')
			+
			\frac{c'+\max\{-c,0\}}{h(\xi_1)}\\
			&\ge
			2+\varepsilon'+\frac{c'}{h(\xi_1)},
		\end{aligned}
		\]	
	\end{proof}
	
	Hence, for $\beta_1,\beta_2\in\Xi'_1$ such that $h(\beta_1)\le h(\beta_2)$,
	we have
	\begin{equation}
		\int_{S\setminus T}
		\min_{1\le i'\le 2}
		\frac{\lambda_{\beta_{i'},J(\omega)}(\omega)}{h(\beta_{i'})}
		\,d\mu(\omega)
		\ge
		2+\varepsilon'
		+\frac{c'}{h(\beta_1)}.
	\end{equation}
	
	Let $$\theta(\omega)=\min_{1\le i'\le 2}
	\frac{\lambda_{\beta_{i'},J(\omega)}(\omega)}{h(\beta_{i'})}.$$
	Then, we have
	\begin{align*}
		\int_{S\setminus T}
		\min_{1\le i'\le 2}
		\theta(\omega)
		\,d\mu(\omega)
		&\ge 2+\varepsilon'
		+\frac{c'}{h(\beta_1)}.
	\end{align*}
	
	We write $\|x\|_\nu=\min\{1,|x|_\nu\}$, so $\|\cdot\|_{\nu}$ is not a norm but still a metric.
	
	Let $j=J(\omega)$. In every $S_{j}$, we have,
	if $\omega$ is Archimedean,
	\begin{align*}
		\|\beta_1-\beta_2\|^{1/\epsilon(\omega)}_\omega
		\le \|\beta_1-\alpha_j\|^{1/\epsilon(\omega)}_\omega
		+\|\beta_2-\alpha_j\|^{1/\epsilon(\omega)}_\omega
		\le 2H(\beta_1)^{-\theta(\omega)/\epsilon(\omega)}.
	\end{align*}
	Next suppose that $\omega$ is non-Archimedean. Then
	$$
	\|\beta_1-\beta_2\|_\omega\
	\le \max\{\|\beta_1-\alpha_j\|_\omega,\|\beta_2-\alpha_j\|_\omega\}
	\le H(\beta_1)^{-\theta(\omega)}.
	$$
	similar to the proof of \ref{amb_gap}, we have 
	$$
	\int_{S\setminus T}-\log^-|\beta_1-\beta_2|_\omega
	\ge (h(\beta_1))\int_{S\setminus T}\theta\,d\mu-h(2)
	\ge (2+\varepsilon'+\frac{c'}{h(\beta_1)})h(\beta_1)-h(2).
	$$
	By Proposition \ref{prop:liou},
	$$
	h(\beta_1)+h(\beta_2)+\log 2
	\ge (2+\varepsilon')h(\beta_1)+c'-h(2).
	$$
	Rearranging the inequality yields
	$$
	h(\beta_2)\ge (1+\varepsilon')h(\beta_1)+c'-h(2)-\log 2.
	$$
	Since $h(2)\le \log 2$, the conclusion follows.
	
	\section{The counting part}\label{counting part}
	Let \(\mathcal G\) be the finite grid constructed in Section \ref{sec:constructive_proof}, consisting of \(M\)
	boxes. For a box \(G\in\mathcal G\), put
	\[
	\mathcal{B}_G:=\mathcal B\cap G,
	\]
	and write its elements in increasing order of height:
	\[
	h(\beta_1)<h(\beta_2)<\cdots<h(\beta_s),
	\qquad s:=\#\mathcal B_G.
	\]
	Choose
	\(\beta_1,\dots,\beta_N\in\mathcal B_G\) with
	$$
	\frac{h(\beta_i)}{h(\beta_{i-1})}\ge r
	\qquad (i=2,\dots,N),
	$$
	and with \(\beta_1\) of minimal height among them. 
	Since \(\beta_1,\dots,\beta_s\) are elements of \(\mathcal B_G\), another version of the first forward gap principle is still satisfied for \(\mathcal{B}_G\) (see Appendix \ref{AppendixA}, which is not so obvious).
	It gives, for every \(j=1,\dots,s-1\),
	\[
	h(\beta_{j+1})\ge(1+\varepsilon')h(\beta_j).
	\]
	
	Set
	\[
	t_0:=\left\lceil\frac{\log r}{\log(1+\varepsilon')}\right\rceil.
	\]
	We claim that
	\begin{equation}\label{number}
		s\le t_0(N-1). 
	\end{equation}
	Suppose in contradiction that
	\[
	s\ge t_0(N-1)+1.
	\]
	Split the sequence \(\beta_1,\dots,\beta_s\) into \(t_0\) subsequences
	according to the residue class of the index modulo \(t_0\):
	\[
	\{\beta_1,\beta_{1+t_0},\beta_{1+2t_0},\dots\},\
	\{\beta_2,\beta_{2+t_0},\dots\},\
	\dots,\
	\{\beta_{t_0},\beta_{2t_0},\dots\}.
	\]
	In each subsequence, consecutive elements are at least \(t_0\) steps apart
	in the original sequence. Hence, by the first-forward gap principle,
	\[
	\frac{h(\beta_{j+t_0})}{h(\beta_j)}
	\ge(1+\varepsilon')^{t_0}
	\ge r
	\]
	for every admissible \(j\). Thus, each subsequence is \(r\)-separated.
	Since
	\[
	s\ge t_0(N-1)+1,
	\]
	the pigeonhole principle yields a subsequence containing at least
	\[
	\left\lceil\frac{t_0(N-1)+1}{t_0}\right\rceil=N
	\]
	elements. All of these \(N\) elements lie in the same box \(G\) and form an
	\(r\)-separated sequence. This contradicts the local box bound (see our detailed proof in p18 and p19), which states that every box in \(\mathcal G\) contains at most
	\(N-1\) elements of any \(r\)-separated subset of \(\mathcal{B}\). Therefore
	(\ref{number}) holds.
	
	Summing (\ref{number}) over all boxes \(G\in\mathcal G\), we obtain
	\[
	\#\mathcal{B}=\sum_{G\in\mathcal G}\#\mathcal{B}_G
	\le t_0(N-1)M.
	\]
	In particular, \(\mathcal{B}\) is finite, and
	\[
	\#\mathcal{B}\le \left\lceil\frac{\log r}{\log(1+\varepsilon')}\right\rceil(N-1)M
	\le \left(1+\frac{\log\Gamma}{\log(1+\varepsilon')}\right)(N-1)M,
	\]
	where the last inequality uses \(r\le\Gamma\).
	
	\vskip0.2cm
	{\footnotesize \noindent
		{\sc YuXi Li}\\
		University of Macau, People's Republic of China.\\
		\textit{E-mail}: mc55457@um.edu.mo}
	\vskip0.2cm
	{\footnotesize \noindent
		{\sc GuanHeng Zhao}\\
		ShanDong University, People's Republic of China.\\
		\textit{E-mail}: 202618373@mail.sdu.edu.cn}
\end{document}